 \documentclass[preprint,review,11pt]{elsarticle}

\usepackage[latin1]{inputenc}
\usepackage[english]{babel}
\usepackage{bm}
\usepackage{graphicx}
\usepackage[nodots]{numcompress}

\usepackage{amsmath}
\usepackage{dsfont}
\usepackage{amssymb}
\usepackage{graphicx}
\usepackage{tabularx}

\newtheorem{theo}{Theorem}
\newtheorem{prop}[theo]{Proposition}

\newtheorem{rem}[theo]{Remark}

\newcommand{\Real}{\mathbb R}

\newcommand{\ds}{\displaystyle}

\newcommand{\bea}{\begin{eqnarray}}
\newcommand{\eea}{\end{eqnarray}}

\journal{}

\begin{document}

\begin{frontmatter}



%========================================================================
% T I T L E
\title{Some numerical aspects of the conservative PSM scheme
 in a 4D drift-kinetic code}
%========================================================================

%========================================================================
% A U T H O R S  and  I N S T I T U T I O N S

\author[INRIA,DIF,IRMA]{Jean-Philippe Braeunig\corref{cor1}}
\ead{braeunig@unistra.fr}
\author[INRIA,IRMA]{Nicolas Crouseilles}
\author[CAD]{Virginie Grandgirard}
\author[CAD]{Guillaume Latu}
\author[IRMA,INRIA]{Michel Mehrenberger}
\author[IRMA,INRIA]{Eric Sonnendr\"{u}cker}

\address[INRIA]{INRIA Nancy Grand-Est, Equipe CALVI,  615 rue du Jardin Botanique F-54600 Villers-l\`es-Nancy, France}
\address[DIF]{CEA, DAM, DIF, F-91297 Arpajon, France}
\address[CAD]{CEA Cadarache, IRFM, F-13108 St Paul-lez-Durance Cedex, France}
\address[IRMA]{IRMA, Universit\'e de Strasbourg, 7 rue Ren\'e-Descartes, F-67084 Strasbourg Cedex, France}
\cortext[cor1]{Corresponding author}

%========================================================================

%========================================================================
% A B S T R A C T   AND   K E Y W O R D S
\begin{abstract}
The purpose of this work is simulation of magnetised plasmas in the ITER project framework. In this context, kinetic Vlasov-Poisson like models are used to simulate core turbulence in the  tokamak in a toroidal geometry. This leads to heavy simulations because a 6D dimensional problem has to be solved, even if reduced to a 5D in so called gyrokinetic models. Accurate schemes, parallel algorithms need to be designed to bear these simulations. This paper describes the numerical studies to improve robustness of the conservative PSM scheme in the context of its development in the GYSELA code. In this paper, we only consider the 4D drift-kinetic model which is the backbone of the 5D gyrokinetic models and relevant to build a robust and accurate numerical method.
\end{abstract}
%========================================================================

\end{frontmatter}

%\linenumbers

%========================================================================
%   I N T R O D U C T I O N
\section{Introduction} \label{sec:introduction}
%- Goal
%- State of the art
 The ITER device is a tokamak designed to study controlled thermonuclear fusion. Roughly speaking, it is a toroidal vessel containing a magnetised plasma where fusion reactions occur. The plasma is kept out of the vessel walls by a magnetic field which lines have a specific helicoidal geometry. However, turbulence develops in the plasma and leads to thermal transport which decreases the confinement efficiency and thus needs a careful study.  Plasma is constituted of ions and electrons, which motion is induced by the magnetic field. The characteristic mean free path is high, even compared with the vessel size, therefore a kinetic description of particles is required, see {\it Dimits} \cite{dimits}. Then the full 6D Vlasov-Poisson model should be used for both ions and electrons to properly describe the plasma evolution. However, the plasma flow in presence of a strong magnetic field has characteristics that allow some physical assumptions to reduce the model. First, the Larmor radius, i.e. the radius of the cyclotronic motion of particles around magnetic field lines, can be considered as small compared with the tokamak size and the gyration frequency very fast compared to the plasma frequency. Thus this motion can be averaged (gyro-average) becoming the so-called guiding center motion. As a consequence, 6D Vlasov-Poisson model is reduced to a 5D gyrokinetic model by averaging equations in such a way the 6D toroidal coordinate system $(r,\theta,\phi,v_\parallel, v_\perp,\alpha)$ becomes a 5D coordinate system $(r,\theta,\phi,v_\parallel, \mu)$, with $v_\parallel$ the parallel and  $v_\perp$ the perpendicular to the field lines components of the particles velocity, $\alpha$ the angular velocity around the field lines and  $\mu=m~v^2_\perp / 2B$ the magnetic momentum depending  on the velocity norm $| v_\perp |$, on the magnetic field magnitude $B$ and on the particles mass $m$ which is an adiabatic invariant. Moreover, the magnetic field is assumed to be steady and the mass of electrons $m_e$ is very small compared to the mass of ions $m_i$. Thus  the cyclotron frequency $\omega_{i,e}=q_{i,e} ~B / m_{i,e}$ is much faster for electrons than for ions  $\omega_{e} >> \omega_{i} $. Therefore the electrons are assumed to be at equilibrium, i.e. the effect of the electrons cyclotronic motion is neglected  and their distribution is then supposed to be constant in time. The 5D gyrokinetic model then reduces to a Vlasov like equation for ions guiding center motion:
\begin{equation}\label{gyro_vlasov}
\begin{array}{ll}
\dfrac{\partial \bar{f}_{\mu}}{\partial t} +   \nabla_X \cdot \left( \dfrac{dX}{dt} \bar{f}_{\mu} \right) + \partial_{v_{\parallel}}  \left( \dfrac{dv_{\parallel}}{dt} \bar{f}_{\mu} \right) = 0
\end{array}
\end{equation}
where $\bar{f}_{\mu} (X,v_\parallel)$ is the ion distribution function for a given adiabatic invariant $\mu$ with $X=(r,\theta,\phi)$, velocities $dX/dt$ and ${dv_{\parallel}}/{dt}$ define the guiding center trajectories. If $\nabla_{(X,v_{\parallel})} \cdot (dX/dt,{dv_{\parallel}}/{dt})^t =0$, then the model is termed as conservative and is equivalent to a Vlasov equation in its advective form:
\begin{equation}\label{gyro_vlasov_nc}
\begin{array}{ll}
\dfrac{\partial \bar{f}_{\mu}}{\partial t} +    \dfrac{dX}{dt} \cdot \nabla_X \bar{f}_{\mu}  +  \dfrac{dv_{\parallel}}{dt} \partial_{v_{\parallel}}  \left( \bar{f}_{\mu} \right) = 0.
\end{array}
\end{equation}
This equation for ions is coupled with a quasi-neutrality equation for the electric potential $\Phi(R)$ on real particles position, with $R=X-\rho_L$ (with $\rho_L$ the Larmor radius) :
\begin{equation}\label{gyro_poisson}
\begin{array}{ll}
-\dfrac{1}{B \omega_{i} } \nabla_\perp \cdot (n_0 \nabla_\perp \Phi) +\dfrac{e}{\kappa T_e}(\Phi-<\Phi>_{\theta,\phi})= \int \bar{f}_{\mu} d\mu dv_{\parallel}-n_0
\end{array}
\end{equation}
where $n_0$ is an equilibrium electronic density, $T_e$ the electronic temperature, $e$ the electronic charge, $\kappa$ the Boltzmann constant for electrons and $ \omega_{i} $ the cyclotronic frequency for ions.\\
These equations are of a simple form, but they have to be solved very efficiently because of the 5D space and the large characteristic time scales considered. This work is then a contribution in this direction, following {\it Grandgirard et al} who develops the GYSELA 5D code that solves this $5D$ gyrokinetic model, see \cite{gysela} and \cite{gysela2}.
Looking at the model, one notices that the adiabatic invariant $\mu$ acts as a parameter. Therefore for each $\mu$ we have to solve a  4D advection equation as accurately as possible but also taking special care on mass and energy conservation, especially in this context of large characteristic time scales. The maximum principle that exists at the continuous level for the Vlasov equation should also be carefully studied at discrete level: $$\ds \min_i ( f(x_i,t^n) ) \le f(x,t^{n+1}) \le \max_i ( f(x_i,t^n) )$$
with  $f(x_i,t^n)$ the value at $x_i$ in cell $i$ at time $t^n$. There is no physical dissipation process in the gyrokinetic model \eqref{gyro_vlasov} that might dissipate over/undershoots created by the scheme and the loss of this bounding extrema of the solution at $t^{n+1} > t^n$ may even eventually crash a simulation.
Those studies will be achieved in this paper on a relevant reduced model, the 4D drift-kinetic model described in section \ref{use},  which has the same structure than equations \eqref{gyro_vlasov}.
The geometrical assumptions of this model for ion plasma turbulence are a cylindrical geometry with coordinates $(r,\theta,z,v_{\|})$ and a constant magnetic field $B=B_z ~e_z$,  where $e_z$ is the unit vector in $z$ direction. %In this collisionless plasma, the trajectories are governed by the guiding center (GC) trajectories:
This 4D model is conservative and will be discretized using a conservative semi-Lagrangian scheme, the Parabolic Spline Method  scheme (PSM, see {\it Zerroukat et al} \cite{PSM1} and \cite{PSM2} ). It is a fourth order scheme which is equivalent for linear advections to the Backward Semi-Lagrangian scheme (BSL) currently used in the GYSELA code (see {\it Grandgirard et al} \cite{gysela2})  and introduced by {\it Cheng-Knorr} \cite{BSL2} and {\it Sonnendr\"{u}cker et al} \cite{BSL}). This conservative PSM scheme based on the conservative form of the Vlasov equation will be described  in section \ref{use} and properly allows a directional splitting. \\
In this paper, the BSL and PSM schemes  will be detailed with an emphasis on their similarities and differences. We will see that one difference is about the maximum principle.  The BSL scheme satisfies it only with a condition on the distribution function reconstruction and  the conservative PSM scheme does not satisfy it  without an extra condition on the volumes conservation in the phase space.
The last condition is equivalent  to try to impose that the velocity field is divergence free at the discrete level. A scheme is given to satisfy this constraint in the form of an equivalent Finite Volume scheme.
Moreover, we have designed a slope limiting procedure, Slope Limited Splines  (SLS), to get closer to a maximum principle for the discrete solution, by at least diminish the spurious oscillations appearing when strong gradients exist in the distribution function profile. \\
The outline of this paper is the following : in section \ref{section_BSLvsPSM} will be recalled some important properties of Vlasov equations at the continuous level. Then BSL and PSM schemes will be described and compared, according to properties of the discrete solutions. In section \ref{stability}, a numerical method will be given to improve the respect of the maximum principle by Vlasov discrete solutions  when using the PSM scheme and particularly to keep constant the volume in the phase space.
%In the same way, a flux limiter will be proposed to improve simulations robustness.
In section \ref{use}, practical aspects of the PSM scheme use will be described  in the context of the 4D drift-kinetic model and at last we will comment on numerical results.
%========================================================================

%#####################################################
\section{Semi-Lagrangian schemes for Vlasov equation}\label{section_BSLvsPSM}
%#####################################################
\subsection{Basics of the Vlasov equation}\label{subsec:basics}
Let us consider an advection equation of a positive scalar function $f(x,t)$ with an arbitrary divergence free velocity field:
\begin{equation}\label{}
\left\{
\begin{array}{l}
\partial_t  f+ a \cdot \nabla_x  (  f )= 0,\\
\nabla \cdot   a =0 ~ \mbox{and} ~f(x,t) \ge 0
\end{array}
\right.
\end{equation}
with position $x \in \Real^D$ and $a(x,t)  \in \Real^D$ the advection velocity field.\\
The solutions satisfy the maximum principle:
\begin{equation}\label{maxprinc}
0 \le f(x,t) \le \ds \max_x (f(x,t_0))
\end{equation}
for any initial time $t_0 < t $.\\
Since $\nabla \cdot   a =0$, we can also  use an equivalent conservative formulation of the Vlasov equation:
\begin{equation}\label{vlasov_consD}
\partial_t  f+ \nabla_x \cdot (  a ~ f )= 0,
\end{equation}
For more details, see {\it Sonnendr\"{u}cker} Lecture Notes \cite{vlasov}.
One obvious property of this conservation law (Reynolds transport theorem) is to conserve  the mass  in a Lagrangian volume $Vol(t)$, by integrating the distribution function on each Lagrangian volume element  $d\Omega$:
 $$d_t m=d_t \int_{Vol(t) }f(x,t) d\Omega=0.$$
Let us introduce the convective derivative $d_t (.)= \partial_t (.)+a \cdot \nabla_x (.)$, thus \eqref{vlasov_consD} becomes:
\begin{equation}\label{}
d_t f + f ~ \nabla_x \cdot  a = 0.
\end{equation}
Considering a Lagrangian motion of an infinitely small volume $Vol(t)$, we have $d_t m = d_t (f ~ Vol) =0$, thus we obtain:
\begin{equation}\label{eqvol}
\dfrac{d_t Vol}{Vol} =  \nabla_x \cdot  a.
\end{equation}
%In particular:
%\begin{equation}\label{eqvol}
%\dfrac{d_t Vol}{Vol} = 0 \Leftrightarrow  \nabla_x \cdot  a=0.
%\end{equation}
Obviously, a divergence free flow $\nabla_x \cdot  a=0 $ conserves a Lagrangian volume in its motion.

%#####################################################
\subsection{Maximum principle in the BSL and PSM schemes}
%=============================================================================
\subsubsection{Backward semi-Lagrangian (BSL)}\label{sec:BSLscheme}
%=============================================================================
Let us consider a Vlasov equation in its non conservative form:
\begin{equation}\label{vlasov_notcons}
\partial_t f+ a \cdot \nabla_x f = 0,
\end{equation}
with $f(x,t)$ a scalar function,  position $x \in \Real^D$ and $a(x,t)  \in \Real^D$ the advection field.
The BSL scheme, see {\it Sonnendr\"{u}cker et al} \cite{BSL}, is based on the invariance property of function $f$ along characteristic curves to obtain values $f^{n+1}$ at time $t^{n+1}$  from the values  $f^n$ at $t^n$:
\begin{equation}
\begin{array}{l}
f^{n+1} \left( X(x^{n+1},t^{n+1}) \right)=f^n \left( X(x^{n+1},t^n) \right),
\end{array}
\end{equation}
with $x$ the Eulerian coordinates and the characteristic curves $X$ defined as
\begin{equation}
\dfrac{d X(x,t) }{dt}=a(x, t)
\end{equation}
with the initial position $x=X(x,t^n)$ at $t^n$.
Let us locate the discrete function values $f_i^n=f^{n} \left( X(x_{i}^{n+1},t^{n+1}) \right)$ at mesh nodes $x_{i}^{n+1}=X(x_{i}^{n+1},t^{n+1}) $. We solve the following nonlinear system which is a second order approximation of $d_t X(t)=a(x,t)$:
\begin{equation}\label{BSLscheme}
\begin{array}{cl}
X_{i}^{n+1/2} & =  \left( X(x_{i}^{n+1},t^{n+1})+X(x_{i}^{n+1},t^n) \right) /2,\\
X(x_{i}^{n+1},t^{n})  &= X(x_{i}^{n+1},t^{n+1})  -  \Delta t ~ a\left( X_{i}^{n+1/2},t^{n+1/2} \right), \\
f^{n+1} \left( X(x_{i}^{n+1},t^{n+1}) \right) & =  f_h^n \left( X(x_{i}^{n+1},t^n) \right),
\end{array}
\end{equation}
with $\Delta t =t^{n+1}-t^n$.
%The scheme consists of two steps:
%\begin{itemize}
%\item At a mesh node $i$ which location is settled as $x_i^{n+1}=X(x_i^{n+1},t^{n+1})$, one has to follow backward the characteristic curve to find the "foot" $x_i^{n}=X(x_i^{n+1},t^{n})$.
%\item
The function $f_h^n(x)$ is a reconstruction of the solution $f^n(x)$  according known values at nodes $x_{i}^{n+1}$ using cubic splines basis functions on the domain to obtain the value  at $x_i^n=X(x_i^{n+1},t^{n})$, which is not a mesh node in general.
%Thus we have $f^{n+1} \left( x_i^{n+1} \right)=f^n \left( x_i^n \right)$.
%\end{itemize}

%=============================================================================
\paragraph{Properties of the BSL scheme}\label{BSL_prop}
%=============================================================================
This scheme is formally fourth order in space. It is second order in time using for instance a Leap-Frog, Predictor-Corrector or Runge-Kutta time integration. Mass is not conserved by this scheme, because it has no conservative form. However, an approximated maximum principle is satisfied. Let us consider $f^n_h(x)$ the cubic spline interpolation of the distribution function $f(x,t^n)$ at time $t^n$, we have for any $x$:
$$f^{n+1} \left( X(x,t^{n+1}) \right) =  f^n_h \left( X(x,t^n) \right).$$
It then naturally appears a "discrete" maximum principle:
\begin{equation}\label{}
\ds \min_x ( f^n_h(x,t^n) ) \le f^{n+1}(x) \le \max_x ( f^n_h(x,t^n) ).
\end{equation}
Comparing with the property \eqref{maxprinc},  we have here $\ds \min_x ( f^n_h(x,t^n) ) \ne 0$ and $\ds \max_x ( f^n_h(x,t^n) )\ne \max_x ( f^n(x,t^n) )$ because the cubic spline reconstruction does not satisfy a maximum principle. If we have a manner to enforce this property to this reconstruction,  a maximum principle is granted for the BSL scheme. No directional splitting is allowed since the BSL scheme is based on the non conservative form of the Vlasov equation, see \cite{vlasov}.
%Directional splitting is not allowed, because considering a $2D$ advection field $(a_x,a_y)$:
%\begin{equation}\label{nospli1}
%\partial_t f + a_x ~\partial_x ( f) +  a_y ~\partial_y ( f) =0
%\end{equation}
%is not equivalent to compute the solution by directional splitting:
%\begin{equation}\label{spli1}
%\left\{
%\begin{array}{ccc}
%\partial_t f + a_x ~\partial_x ( f)  =0 \\
%\partial_t f + a_y ~\partial_y ( f) =0.
%\end{array}
%\right.
%\end{equation}
%Even if the velocity  field $(a_x,a_y)^t$ is such that $\nabla \cdot (a_x,a_y)^t=0$, each step of (\ref{spli1}) is not conservative because $\partial_x a_x \ne 0$ and  $\partial_y a_y \ne 0$ in general, see \cite{BSLinst} to go further. A directional splitting is only allowed when using the conservative form (\ref{vlasov_cons}), see section $\ref{section_PSM}$.

%=============================================================================
\subsubsection{Semi-Lagrangian Parabolic Spline Method (PSM)}\label{section_PSM}
%=============================================================================
Let us consider a Vlasov equation in its conservative form:
\begin{equation}\label{vlasov_cons}
\partial_t f + \nabla_x \cdot (  a ~ f )= 0,
\end{equation}
with $f(x,t)$ a scalar function,  position $x \in \Real^D$ and $a(x,t)  \in \Real^D$ the advection field. Notice that with the hypothesis $\nabla_x \cdot (  a )=0$, conservative form \eqref{vlasov_cons} and non-conservative  form \eqref{vlasov_notcons} of the Vlasov equation are equivalent.
The PSM scheme, see {\it Zerroukat et al} \cite{PSM1} and  \cite{PSM2}, is based on the mass conservation property of function $f$ in a Lagrangian  volume $Vol$ to obtain the value $f^{n+1}$ at time $t^{n+1}$:
\begin{equation}\label{SLPSM}
\begin{array}{l}
\displaystyle \int_{Vol^{n+1}} f(x,t^{n+1}) d\Omega=\displaystyle \int_{Vol^{n}} f(x,t^{n}) d\Omega,
\end{array}
\end{equation}
with the characteristic curves $X$ defined as $\dfrac{d X(x,t) }{dt}=a(x, t)$ and $x^n=X(x,t^n)$ with $t^n$ the initial time, and the volume
$Vol^n=\{ X(x^{n+1},t^n)$ such that $X(x^{n+1},t^{n+1}) \in Vol^{n+1} \}$ defined by the Lagrangian motion with the field $a(x,t)$.
The important point is that this conservative formalism properly allows a directional splitting without loosing the mass conservation. Indeed, equation \eqref{vlasov_cons} may be solved with $D$ successive 1D advections still of conservative form:
\begin{equation}\label{}
\partial_t f + \partial_{x_k} (  a_k ~ f )= 0,~ k \in [1,D].
\end{equation}
We then approximate a 1D equation for each direction $k$ using the conservation property. Omitting subscript $k$, the PSM scheme writes in 1D as follows:
\begin{equation}\label{cons_1D}
\begin{array}{l}
\displaystyle \int^{ x_{i+1/2}^{n+1} }_{x_{i-1/2}^{n+1} } f(x,t^{n+1}) dx=\displaystyle \int_{x_{i-1/2}^{n}}^{x_{i+1/2}^{n}} f(x,t^{n}) dx,
\end{array}
\end{equation}
with $x_{i+1/2}^{n+1}=X(x_{i+1/2}^{n+1},t^{n+1}) $ settled as the 1D mesh nodes and $x_{i+1/2}^{n}=X(x_{i+1/2}^{n+1},t^{n})$ the associated foot of the characteristic curve, $Vol_i^n = [x^n_{i-1/2}, x^n_{i+1/2}] $ and  $Vol_i^{n+1} = [x^{n+1} _{i-1/2}, x^{n+1} _{i+1/2}]$.\\
Let us define the unknowns of the scheme as the average of $f$ in cell $i$
\begin{equation}
\begin{array}{l}
\overline{f}^{n+1}_i=\dfrac{1}{\Delta x} \displaystyle \int^{ x_{i+1/2}^{n+1} }_{x_{i-1/2}^{n+1} } f(x,t^{n+1}) dx,
\end{array}
\end{equation}
 and the primitive function
\begin{equation}
\begin{array}{l}
F^{n} (z)=\displaystyle \int^{ z }_{ x_{1/2}} f(y,t^{n}) dy,
\end{array}
\end{equation}
with the uniform space step $\Delta x=x^{n+1} _{i+1/2}-x^{n+1} _{i-1/2}$ and $x_{1/2}$ an arbitrary reference point of the domain and for instance the first node of the grid $\{x_{i-1/2}\}_{i=1,N+1}$.
Therefore, one has to solve a nonlinear system, which is similar to the BSL one, to obtain a discrete solution of equation \eqref{cons_1D} that writes:
\begin{equation}\label{PSMscheme}
\begin{array}{cl}
X_{i+1/2}^{n+1/2} & =  \left( X(x_{i+1/2}^{n+1},t^{n+1})+X(x_{i+1/2}^{n+1},t^n) \right) /2,\\
X(x_{i+1/2}^{n+1},t^{n})  &= X(x_{i+1/2}^{n+1},t^{n+1})  -  \Delta t ~ a\left( X_{i+1/2}^{n+1/2},t^{n+1/2} \right), \\
\overline{f}^{n+1}_i \Delta x &=  F_h^{n} (X(x_{i+1/2}^{n+1},t^{n}))-F_h^{n}(X(x_{i-1/2}^{n+1},t^{n}))
\end{array}
\end{equation}
with the time step $\Delta t =t^{n+1}-t^n$ and the uniform space step $\Delta x=x_{i+1/2}^{n+1}-x_{i-1/2}^{n+1}$ . \\
%The scheme is constituted with two steps:
%\begin{itemize}
%\item  For a mesh node $i+1/2$ which location is settled as $x_{i+1/2}^{n+1}=X(x_{i+1/2}^{n+1},t^{n+1})$, one has to follow backward the characteristic curve to find the "foot" $x_{i+1/2}^{n}=X(x_{i+1/2}^{n+1},t^{n})$.
%\item
The computation of the reconstructed primitive function $F^n_h(x)$  is based on values at mesh nodes $x_{i+1/2}^{n+1}$:
$$F_h^{n} ( x_{i+1/2}^{n+1} ) - F_h^{n} ( x_{1/2} )=\sum^{i}_{k=1} \overline{f}^{n}_k \Delta x.$$
Then this set of values is interpolated by cubic splines functions to obtain an approximated value $F_h^{n} (z)$  of the primitive function $F^{n} (z)$  at any point $z$ of the domain:
\begin{equation}\label{PSMprimitive}
\begin{array}{cl}
F_h^{n} (z) \approx F^{n} (z) = \displaystyle \int^{ z }_{ x_{1/2}} f(y,t^{n}) dy.
\end{array}
\end{equation}

%=============================================================================
\paragraph{Properties of  the PSM scheme}
%=============================================================================
This scheme is formally fourth order in space and strictly equivalent to the BSL scheme for constant linear advection, see \cite{PSM}. It is second order in time using for instance a Leap-Frog, Predictor-Corrector or Runge-Kutta time integration scheme. Mass is exactly conserved by this scheme for each  1D step $k$ of the directional splitting.
%\begin{equation}\label{conservation}
%\begin{array}{ll}
%\displaystyle \int^{ x_{N+1/2}^{n+1} }_{x_{1/2}^{n+1} } f(x,t^{n+1}) dx
%& =  \ds \sum^{N}_{k=1} \overline{f}^{n+1}_k \Delta x =  F^{n} (x_{N+1/2}^n)-F^{n} (x_{1/2}^n)\\
%~ & = \displaystyle \int_{x_{1/2}^{n}}^{x_{N+1/2}^{n}} f(x,t^{n}) dx.
%\end{array}
%\end{equation}
However, no maximum principle does exist  for each step $k$ even for the exact solution:  in general $\partial_{x_k} a_k \ne 0$, even if the velocity field is divergence free $\nabla \cdot a=0$.
%, so the maximum principle for non constant advections only exists considering all steps of the dimensional splitting.
Let us consider the scheme in $D$ dimensions of space:
\begin{equation}\label{vlas3}
\partial_t f + \nabla \cdot (  a ~ f )= 0,
\end{equation}
with $f(x,t)$ a scalar function,  position $x \in \Real^D$ and $a(x,t)  \in \Real^D$ the advection field. Let us consider a cell $i$, where the solution is described at time $t^{n+1}$ by its average in cell $i$,
%$$\overline{f}^{n+1}_i =\dfrac{1}{Vol_i^{n+1}} \displaystyle \int_{Vol_i^{n+1}} f(x,t^{n+1}) d\Omega,$$
%with $Vol_i^{n+1}$ settled as the cell $i$ volume.\\
 the PSM scheme then writes:
\begin{equation}
\begin{array}{l}
\overline{f}^{n+1}_i Vol_i^{n+1}=\displaystyle \int_{Vol_i^{n+1}} f(x,t^{n+1}) d\Omega=\displaystyle \int_{Vol_i^{n}} f(x,t^{n}) d\Omega.
\end{array}
\end{equation}
%with the characteristic curves $X$ defined as $\dfrac{d X(x,t) }{dt}=a(x, t)$ and $x=X(x,t^n)$ with $t^n$ the initial time, and the volume
%$Vol_i^n=\{ X(x^{n+1},t^n)$ such that $X(x^{n+1},t^{n+1}) \in Vol_i^{n+1} \}$ defined by the Lagrangian motion with the field $a(x,t)$.
We thus obtain the following relation:
\begin{equation}
\overline{f}^{n+1}_i =\overline{f}^{n*}_i \dfrac{Vol_i^{n}}{ Vol_i^{n+1}}
\end{equation}
with the average of the distribution function in the Lagrangian volume at time $t^n$:
$$ \overline{f}^{n*}_i =\dfrac{1}{Vol_i^{n}} \displaystyle \int_{Vol_i^{n}} f(x,t^{n}) d\Omega.$$

Here clearly appears two conditions, both difficult to satisfy especially in the context of a directional splitting, to have a maximum principle defined as follows:
\begin{equation}\label{max_princ}
\min_j ( \overline{f}^{n}_j ) \le \overline{f}^{n+1}_i \le \max_j ( \overline{f}^{n}_j ).
\end{equation}

\begin{enumerate}
\item Maximum principle on the distribution function in $Vol_i^{n}$:
\begin{equation}\label{mono_aver}
\min_j ( \overline{f}^{n}_j ) \le \overline{f}^{n*}_i \le \max_j ( \overline{f}^{n}_j ).
\end{equation}
\item Conservation of volumes in the phase space at the discrete level:
\begin{equation}\label{cons_vol}
Vol_i^{n} = Vol_i^{n+1}.
\end{equation}
\end{enumerate}
The first condition is difficult to ensure in general, because a maximum principle should be satisfied for any average of the distribution function on an arbitrary volume $Vol_i^n$. Moreover, in the context of a directional splitting, it is impossible to satisfy a maximum principle for a 1D step $k$, because it does not exist at the continuous level since in general $\partial_k a_k \ne 0$. Therefore it is probably impossible to recover a maximum principle of the reconstruction after all steps of the directional splitting.\\
The second condition is true at the continuous level while $\nabla \cdot a=0$, since we have $d_t Vol = Vol ~ \nabla \cdot a=0$, see equation \eqref{eqvol} in section \ref{subsec:basics}. As well as for the first condition, in the context of a directional splitting, it is difficult to ensure a constant volume evolution $Vol_i^{n} = Vol_i^{n+1}$ after all steps of the directional splitting, where compressions or expansions of the Lagrangian volume occur successively.  \\

As a consequence, we will propose a  form of the conservative PSM scheme that does not use a directional splitting. However we will not write the  Semi-Lagrangian form of the  PSM scheme in $D$ dimensions of space because it is costly in computational time, because of the reconstruction step, and it is difficult to handle with arbitrary coordinate systems. The solution we choose is to use an equivalent Finite Volume form of the PSM scheme described in section \ref{PSM_FV}, which is locally 1D at each face of the mesh. It is therefore possible to design 1D numerical limiters to try to better satisfy the maximum principle condition \eqref{mono_aver}. Moreover, we will show that this form allows an exact conservation of the volumes in the phase space \eqref{cons_vol}.  The maximum principle and therefore the robustness of this scheme will thus be considerably improved.

%#####################################################%#####################################################
\section{Maximum principle for the PSM scheme}\label{stability}
%#####################################################%#####################################################

%#####################################################
\subsection{Numerical limiters for the distribution function reconstruction}\label{consvol}
\begin{figure}%[!h]
 \centering
 \begin{tabular}{c}
   \includegraphics[scale=0.35]{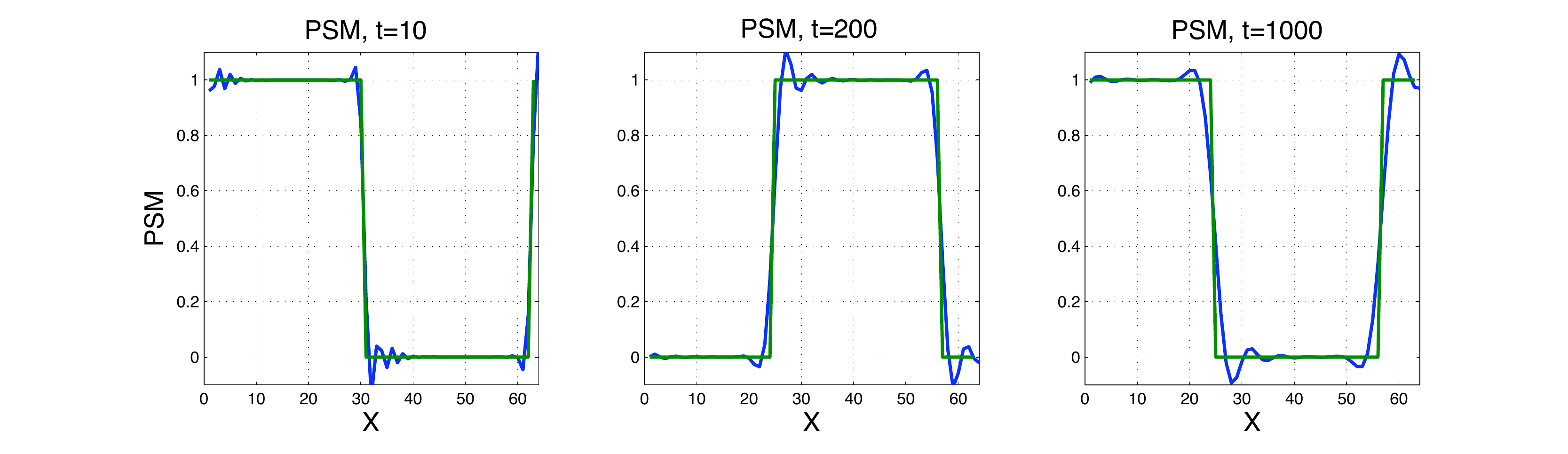} \\
   \includegraphics[scale=0.35]{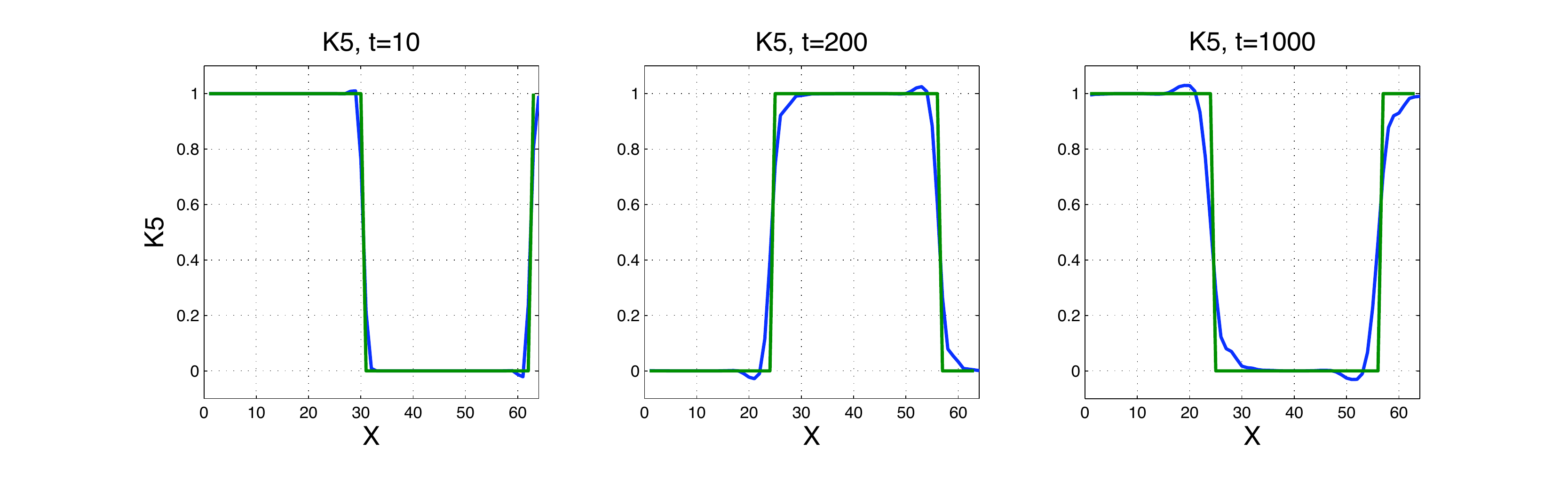}
 \end{tabular}
 \caption{ \label{fig:advstep}Linear advection of a step function. The exact solution is green and the numerical results are in blue, above the standard PSM scheme and below the PSM scheme with SLS limiter with $K=5$. }
 \end{figure}

Enforcing the first condition on the maximum principle of the distribution function reconstruction \eqref{mono_aver} can be really costly in computational time. Instead of trying to correct the cubic spline reconstruction, we will reduce the spurious oscillations, generated by high order schemes when strong gradients appear in the distribution function profile, by using a classical {\it Van Leer} like slope limiting procedure, see for instance {\it LeVeque} \cite{leveque}.
We propose here to measure the gradients in the flow and to add diffusion where they are detected. The diffusion is added by mixing the high order PSM flux with a first order upwind flux. The evaluation of the gradient is given by the classical function $\theta$ and we estimate the diffusion needed with a function $\gamma(\theta) \in [0,1]$ based on a minmod like limiter function (see Fig. \ref{fig:gamma_theta}). The resulting limiter we propose here is called SLS (Slope Limited Splines), see \cite{guterl} for details:
\[
\phi^{SLS}_{i+1/2}=\gamma(\theta_{i+1/2}) ~\phi^{PSM}_{i+1/2}+(1-\gamma(\theta_{i+1/2})) ~\phi^{upwind}_{i+1/2}
\]
 where
$$
\phi^{upwind}_{i+1/2}= a _{i+1/2} \left( \dfrac{\overline{f}^{n}_{i}+\overline{f}^{n}_{i+1}}{2}- \mbox{sign}(a_{i+1/2}) \dfrac{\overline{f}^{n}_{i+1}-\overline{f}^{n}_{i}}{2} \right).
$$
We define $\theta_{i+1/2}$ as the classical slope ratio of the distribution which depends on the direction of the displacement:
\[ \theta_{i+1/2}= \left\{
    \begin{array}{ll}
        \dfrac{\overline{f}^{n}_{i}-\overline{f}^{n}_{i-1}}{\overline{f}^{n}_{i+1}-\overline{f}^{n}_{i}} & \text{ if } a_{i+1/2}>0  \\
      \dfrac{\overline{f}^{n}_{i+2}-\overline{f}^{n}_{i+1}}{\overline{f}^{n}_{i+1}-\overline{f}^{n}_{i}} & \text{ if } a_{i+1/2}<0   \end{array}
\right.
\]
%check1
%\begin{figure}%[!h]
However, the classical limiter minmod, where $\gamma_{i+1/2}=\max(0,\min(\theta_{i+1/2},1))$, set $\gamma$ to 0 when $\theta<0$. That means that the scheme turns to order 1 when an extrema exists, i.e. the slope ratio $\theta<0$. These extrema are thus quickly diffused and that leads to loose the benefits of a high order method.
For SLS, the choice is to let the high-order scheme deal with the extrema and only add diffusion when strong gradients occurs, i.e.  the slope ratio $\theta \approx 0$. We also introduce a constant K in relation to control the maximum slope allowed without adding diffusion, i.e. mixing with the upwind scheme, see figure \ref{fig:gamma_theta}:
 \begin{figure}%[!h]
 \centering
 \begin{tabular}{c}
   \includegraphics[scale=0.35]{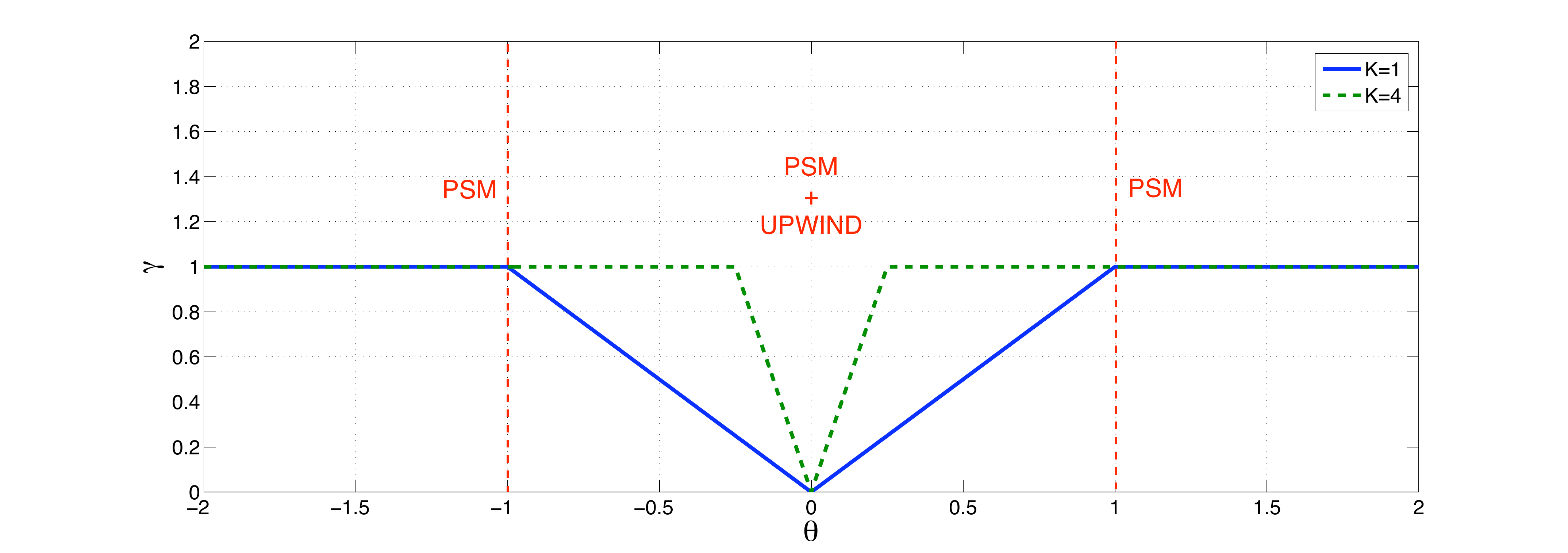} \\
 \end{tabular}
 \caption{$\gamma$ function for the SLS limiter}
 \label{fig:gamma_theta}
 \end{figure}
\begin{equation}\label{sls}
\gamma_{i+1/2}=\max(0,\min(K|\theta_{i+1/2}|,1)).
\end{equation}
with the constant $K=5$ experimentally settled. \\
We present in Fig. \ref{fig:advstep} the results of the linear advection of a step function with the standard PSM scheme with and without the SLS limiter (K=5). The domain is meshed with 70 cells with periodic boundary conditions and the displacement is set to 0.2 cell per iteration.
One can see that as any high order scheme, the PSM scheme produces spurious oscillations at the discontinuity or at a stiff gradient location. The SLS limiter do well with $K=5$  to reduce these oscillations without introducing much diffusion. However, a maximum principle is not granted. This limiter has been further  studied and compared with other limiters in the report \cite{guterl}.

%#####################################################
\subsection{Finite Volume form of the PSM scheme}\label{PSM_FV}

 Let us consider the 1D conservative advection equation of the form:
\begin{equation}\label{vlasov_cons1D}
\partial_t f + \partial_x (  a_x ~ f )= 0,
\end{equation}
with $f(x,t)$ a scalar function,  position $x \in \Real$ and $a_x(x,t)  \in \Real$ the advection field.
We recall that  $X(x_{i+1/2}^{n+1},t^{n+1})=x_{i+1/2}$ is the position of the mesh node $i+1/2$. Let us set the notation $X(x_{i+1/2}^{n+1},t^{n})=x^*_{i+1/2}$ for the "foot" position at $t^n$ on the characteristic curve.  Let us rewrite the PSM scheme \eqref{PSMscheme}:
$$x^*_{i+1/2}= x_{i+1/2} - \Delta t ~ a_x \left( x^{n+1/2}_{i+1/2}  , t^{n+1/2} \right),$$
$$ \overline{f}^{n+1}_i \Delta x = F_h^{n} (x^*_{i+1/2})-F_h^{n}(x^*_{i-1/2}).$$
with the primitive function $F_h^{n} (z)$ is interpolated by cubic splines at mesh nodes such that $F_h^{n} (z)   \approx F^{n} (z) =\int^{z}_{x_{1/2}} f^n(y) dy$ at fourth order in space.\\
Let us make appear explicitly in 1D the fluxes at cell faces $i\pm1/2$, by introducing the primitive values at cell faces $F^{n} (x_{i\pm 1/2})$:
$$ \overline{f}^{n+1}_i \Delta x = \left( F_h^{n} (x^*_{i+1/2}) - F_h^{n} (x_{i+1/2}) \right)
- \left(F_h^{n}(x^*_{i-1/2}) - F_h^{n}(x_{i-1/2}) \right)+\overline{f}^{n}_i \Delta x $$
with $\overline{f}^{n}_i \Delta x=F_h^{n}(x_{i+1/2})-F_h^{n}(x_{i-1/2})$. It yields
\begin{equation}\label{PSM_form_flux}
\dfrac{ \overline{f}^{n+1}_i - \overline{f}^{n}_i } {\Delta t}  +   \dfrac{ F_h^{n} (x_{i+1/2}) - F_h^{n} (x^*_{i+1/2}) }{\Delta x \Delta t}
-    \dfrac{ F_h^{n} (x_{i-1/2}) - F_h^{n} (x^*_{i-1/2}) }{\Delta x \Delta t}  =0.
\end{equation}
The PSM fluxes  at cell faces $i\pm1/2$ clearly appear:
\begin{equation}\label{form_flux}
\dfrac{ \overline{f}^{n+1}_i - \overline{f}^{n}_i } {\Delta t}  +   \dfrac{ \Phi^{PSM}_{i+1/2} - \Phi^{PSM}_{i-1/2} }{\Delta x}   =0
\end{equation}
with
\begin{equation}\label{flux_PSM}
\Phi^{PSM}_{i+1/2} =  \dfrac{ F_h^{n} (x_{i+1/2}) - F_h^{n} (x^*_{i+1/2}) }{ \Delta t} \approx\dfrac{1}{\Delta t} \int^{x_{i+1/2}}_{x^*_{i+1/2}} f^n(y) dy.
\end{equation}
A simple Taylor expansion shows that this PSM flux, which consists in a cubic spline approximation of the integral of $f^n (x)$ along the characteristic curves at cell faces, is a consistent approximation at node $i+1/2$ of the continuous flux $\Phi=a_x f$ in equation  \eqref{vlasov_cons1D}, i.e. $\Phi^{PSM}_{i+1/2}  \approx (a_x f)_{i+1/2}.$\\
Moreover, this flux is an approximation of the integral of $(a_x f^n$ on cell faces. Coming back to \eqref{vlasov_cons1D} and integrating on the cell volume $Vol_i^{n+1}=A_x \Delta x$,  with $A_x$ the bounding faces $\Gamma_{i\pm1/2}$ area transversal to $x$ of cell $i$:
\begin{equation}\label{vlasov_cons1Db}
Vol_i^{n+1} \frac{\partial \overline{f}}{\partial t} + \int_{Vol_i^{n+1}}\partial_x (  a_x ~ f ) d\Omega= 0.
\end{equation}
We obtain using {\it Green} formula:
\begin{equation}\label{vlasov_cons1cd}
\frac{\partial \overline{f}}{\partial t} + \dfrac{1}{Vol_i^{n+1}}\int_{\partial Vol_i^{n+1}=\Gamma_{i-1/2} \cup \Gamma_{i+1/2}} f (  a_x \cdot n_x ) d\Gamma= 0.
\end{equation}
with $\partial Vol_i^{n+1}=\Gamma_{i-1/2} \cup \Gamma_{i+1/2}$ the surface bounding $Vol_i^{n+1}$ and $n_x$ its outgoing normal.
\begin{equation}\label{vlasov_cons1c}
\frac{\partial \overline{f}}{\partial t} + \dfrac{1}{Vol_i^{n+1}} \left( \int_{\Gamma_{i+1/2} } f a_x  d\Gamma-\int_{\Gamma_{i-1/2} } f a_x  d\Gamma \right)= 0.
\end{equation}
Comparing with formula \eqref{form_flux}, we see that the PSM flux is an approximation of the flux average at cell faces:
\begin{equation}\label{flux1}
 \Phi^{PSM}_{i+1/2}\approx\dfrac{1}{A_x}\int_{\Gamma_{i+1/2}} f (  a_x \cdot n_x ) d\Gamma
\end{equation}
because  $Vol_i^{n+1}=A_x \Delta x$.\\

The extension to $D$ dimensions of space is then straightforward, because it only consists in adding the fluxes through the faces of a cell in every direction $d$. Considering  the Vlasov equation in its conservative form in dimension $D$ and for an arbitrary coordinate system:
$$\frac{\partial( J  f)}{\partial t} + \nabla \cdot ( J  a ~ f )= 0$$
with $J$ the geometric Jacobian of the cell, $f(x,t)$ a scalar function,  position $x \in \Real^D$ and $a(x,t)  \in \Real^D$ the advection field. In a classical way in Finite Volume methods, we integrate this local equation on the cell $C_i$ of volume $Vol_i^{n+1}$ and we use the {\it Green} formula:
\begin{equation}
\displaystyle \int_{Vol^{n+1}} \dfrac{\partial f}{\partial t} J dx + \displaystyle \int_{\Gamma_{d} \in C_i}  (a \cdot n_d) f  J dx =0,
\end{equation}
with $\Gamma_{d} $ the face of cell $i$ perpendicular to direction $d$ of area $A_{d}$ and of outgoing normal unit vector $n_d$.
Using the 1D flux formula \eqref{flux1} in direction $n_d$ and a first order time discretisation, we obtain the following Finite Volume scheme:
\begin{equation}\label{form_fluxD}
Vol_i^{n+1} \dfrac{ \overline{f}^{n+1}_i - \overline{f}^{n}_i } {\Delta t}  + \displaystyle \sum_{\Gamma_{d} \in C_i} A_{d}  \Phi^{PSM}_{d}    =0,
\end{equation}
with
\begin{equation}\label{form_fbD}
 \overline{f}^{n}_i=\dfrac{1}{Vol_i^{n+1}} \displaystyle \int_{Vol_i^{n+1}}  f(y,t^{n}) J(y) dy.
\end{equation}
and $\Phi^{PSM}_{d} $ the flux that goes through $\Gamma_{d} $:
\begin{equation}\label{flux_PSMD}
\Phi^{PSM}_{d}  =  \dfrac{ F_{h,d}^{n} (x_{d}) - F_{h,d}^{n} (x^*_{d}) }{ \Delta t} \approx  \dfrac{1}{A_{d}} \displaystyle \int_{\Gamma_{d} }  (a \cdot n_d) f(y,t^n)  J(y)  dy,
\end{equation}
with $x^*_{d} = x_{d} - \Delta t ~ (a(x_d) \cdot n_d)$ the foot of the 1D characteristic curve and $F^n_{h,d}(z)$ the primitive function at $t^n$ reconstructed using cubic splines in the direction of $n_{d}$:
\begin{equation}\label{PSMprimitived}
\begin{array}{cl}
F_{h,d}^{n} (z) \approx F_d^{n} (z) = \displaystyle \int^{ z_d }_{ x_{d,1/2}}  f(y,t^{n}) J(y) dy.
\end{array}
\end{equation}
We see here a Finite Volume form of the Semi-Lagrangian PSM conservative scheme. This equivalence is however restricted, because this Finite Volume form is submitted to a CFL condition as any scheme of this form:
$$\Delta t \le  \displaystyle \min_d \left( \dfrac{\Delta x_d}{ \displaystyle \max_{x_d} (a^n_d(x_d))} \right).$$
Moreover, as we will see in section \ref{consvol}, the Lagrangian volume evolution is here approximated by the cell faces motion only in their normal direction, instead of the general motion as it is described in the Semi-Lagrangian formalism \eqref{SLPSM}. It is the same volume evolution  as  the Semi-Lagrangian method with 1D directional splitting. However it is the classical Finite Volume formalism and it is the key point that will permit to enforce a divergence free evolution of the flow.\\
Notice that the Finite Volume form \eqref{form_fluxD} can be directionally split keeping exactly the same result. Indeed, it only consists in adding the flux in two successive operations instead of in one. As an example, let us consider the 2D $(x,y)$ case of a cartesian mesh:
\begin{equation}\label{form_fluxDxy}
   \begin{array}{ll}
Vol_i^{n+1} \dfrac{ \overline{f}^{n+1}_i - \overline{f}^{n}_i } {\Delta t}
& + \displaystyle  A^x_{i+1/2,j}  \Phi^{PSM}_{i+1/2,j} + A^x_{i-1/2,j}  \Phi^{PSM}_{i-1/2,j}  \\
~ &+A^y_{i,j+1/2}  \Phi^{PSM}_{i,j+1/2} +A^y_{i,j-1/2}  \Phi^{PSM}_{i,j-1/2}   =0.
   \end{array}
\end{equation}
It is stricly equivalent to use the  directional directional  splitting:
\begin{equation}\label{form_fluxDxy_split}
   \begin{array}{ll}
Vol_i^{n+1} \dfrac{ \overline{f}^{nx}_i - \overline{f}^{n}_i } {\Delta t}
& + \displaystyle  A^x_{i+1/2,j}  \Phi^{PSM}_{i+1/2,j} + A^x_{i-1/2,j}  \Phi^{PSM}_{i-1/2,j} =0 \\
Vol_i^{n+1}  \dfrac{ \overline{f}^{n+1}_i - \overline{f}^{nx}_i } {\Delta t}   &+A^y_{i,j+1/2}  \Phi^{PSM}_{i,j+1/2} +A^y_{i,j-1/2}  \Phi^{PSM}_{i,j-1/2}   =0.
   \end{array}
\end{equation}
with the only condition that all fluxes $\Phi^{PSM}$ are computed using $f^n$ and $a^n$ at time $t^n$ as in the unsplit scheme \eqref{form_fluxDxy}.

%#####################################################
\subsection{Conservation of volumes in the phase space}\label{consvol}
\begin{figure}%[!h]
 \centering
 \begin{tabular}{c}
   \includegraphics[scale=0.7]{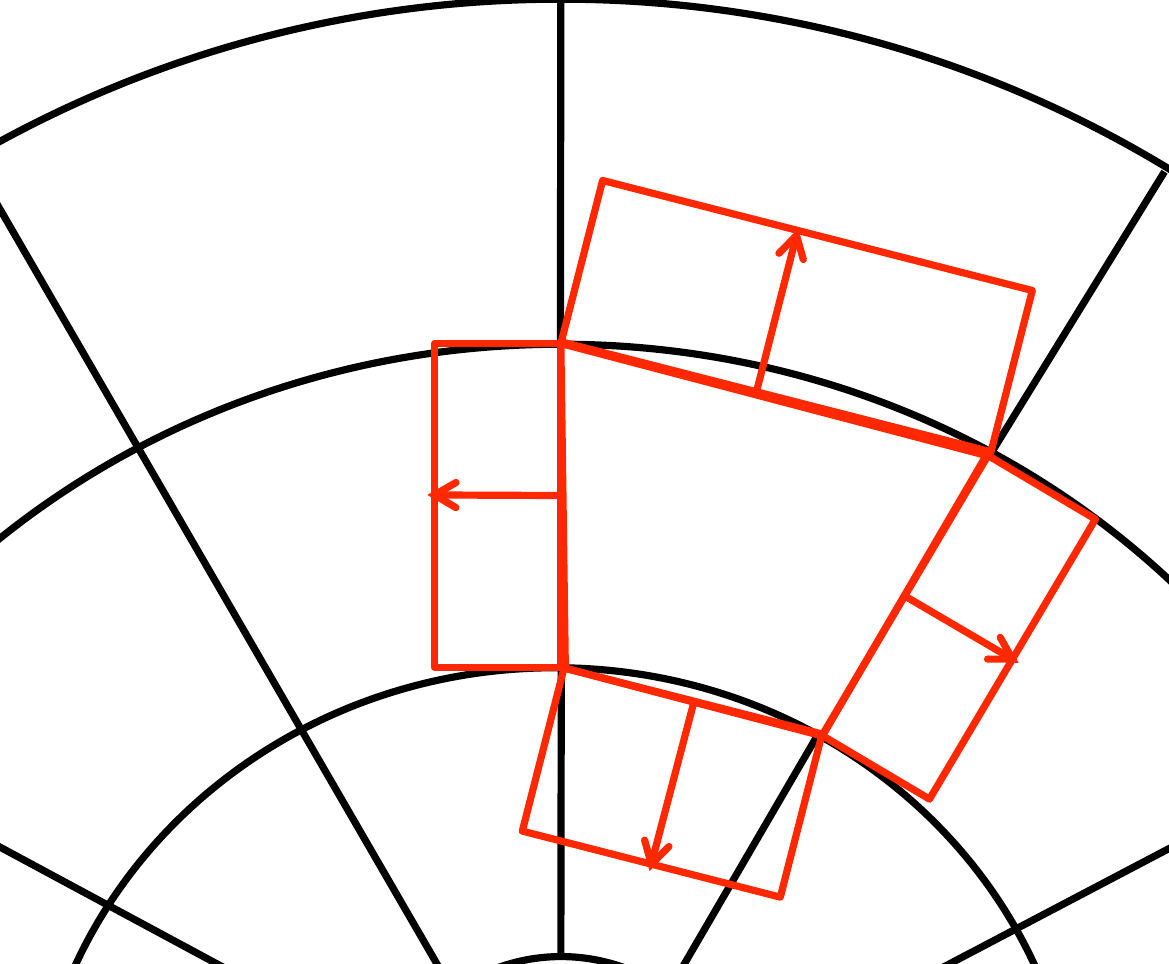} \\
 \end{tabular}
 \caption{Polar cells in black and the Finite Volume representation of the cells in red as a trapezium with the volumes swept by the cell faces in their normal motion. The direction for all faces motion is drawn outward, but it could be inward as well, function of the velocity field.}
 \label{fig:polVF}
 \end{figure}

The second condition to have a maximum principle for the PSM scheme is to satisfy the multi-dimensional condition \eqref{cons_vol} of conservation of volumes, i.e. $Vol^n=Vol^{n+1}$. Equation \eqref{eqvol} showed that  at continuous level the volume is constant in its evolution in the phase space if the advection field is divergence free.  Therefore we will study the PSM scheme to find out a divergence free condition that should be satisfied at the discrete level $\nabla^h \cdot a=0$ in such a way $Vol^n=Vol^{n+1}$, in the same way $\nabla \cdot a=0$ at the continuous level. With the idea of making appear the total evolution of volumes between $t^n$ and $t^{n+1}$, we will use the Finite Volume form of the PSM scheme given in section \ref{PSM_FV} for the 2D polar coordinate system.
 Let us consider radial $r$ and orthoradial $\theta$ directions, with constant space steps $\Delta r$, $\Delta \theta$ and the volume of the cells $Vol_{i,j}=r_{i} \Delta r \Delta \theta$ with $r_i$ the mean radius of the cell. We consider that the mesh in polar coordinates has locally no curvature, i.e. each mesh is a trapezium with straight edges. This is important to be noticed to write the Finite Volume scheme and calculate the volume swept by the cell edges, see Fig. \ref{fig:polVF}.
 Let us set a velocity field $\left( a_r(r,\theta), r ~ a_\theta(r,\theta) \right)$ such that:
$$\nabla_{r,\theta} \cdot a= \dfrac{1}{r} \partial_r (r ~a_r) +  \dfrac{1}{r} \partial_\theta (r ~a_\theta) =0.$$
Let us write the conservative advection equation in polar coordinates:
\begin{equation}\label{}
\partial_t (r f) +  \partial_r (r ~ a_r ~f)  +  \partial_\theta (r ~ a_\theta ~f)=0.
\end{equation}
Notice that the geometric Jacobian $J=r$ for polar coordinates. \\
The PSM scheme without directional splitting in the Finite Volume form \eqref{form_fluxD} reads here:
\begin{equation}\label{ex_split_2Drt}
%\left\{
\begin{array}{ll}
Vol_{i,j} \dfrac{ \overline{f}^{n+1}_{i,j} - \overline{f}^{n}_{i,j} } {\Delta t}      & + A^r_{i+1/2,j} \Phi^{PSM,r}_{i+1/2,j} - A^r_{i-1/2,j}  \Phi^{PSM,r}_{i-1/2,j}    \\
~ & + A^\theta_{i,j+1/2} \Phi^{PSM,\theta}_{i,j+1/2} - A^\theta_{i,j+1/2}  \Phi^{PSM,\theta}_{i,j-1/2}   =0
\end{array}
%\right.
\end{equation}
with $\Phi^{PSM,\theta}_{i,j\pm1/2}$ and $\Phi^{PSM,r}_{i\pm1/2,j}$ positioned at cell faces center and with cell $i,j$ of volume $Vol_{i,j}=r_i \Delta r \Delta \theta$ and faces areas $A^r_{i \pm 1/2,j}=\Delta \theta$ and $A^\theta_{i,j\pm1/2}=\Delta r$.
The cell averaged values of $ \overline{f}$ used in the scheme are:
\begin{equation}
\overline{f}_{i,j} =\dfrac{1}{Vol_{i,j}} \int_{Vol_{i,j}} f(r,\theta,t)~ r dr d\theta.
\end{equation}
Using the integral form \eqref{flux_PSMD} of  the fluxes:
\begin{equation}\label{}
%\left\{
\begin{array}{cl}
~ & r_i \Delta r \Delta \theta ( \overline{f}^{n+1}_{i,j} -  \overline{f}^{n}_{i,j} ) \\
 ~ & + \displaystyle  \Delta \theta \int^{r_{i+1/2}}_{r^*_{i+1/2}} f(r,\theta_i,t^{n}) r dr -\displaystyle  \Delta \theta \int^{r_{i-1/2}}_{r^*_{i-1/2}} f(r,\theta_i,t^{n}) r dr   \\
~ & + \displaystyle \Delta r \int^{\theta_{j+1/2}}_{\theta^*_{j+1/2}} f(r_i,\theta,t^{n}) r_i d\theta - \Delta r\displaystyle \int^{\theta_{j-1/2}}_{\theta^*_{j-1/2}} f(r_i,\theta,t^{n}) r_i  d\theta  =0.
\end{array}
%\right.
\end{equation}
Let us introduce the volumes swept by each cell face in its normal motion in accordance with the {\it Green} formula and the way of computation of feet of characteristic curves normal to cell faces, i.e. without taking into account the tangential motion at the cell faces or the curvature of the mesh, see Fig. \ref{fig:polVF}:
$$\delta Vol^r_{i \pm 1/2,j} =  r_{i\pm1/2}  ~ \Delta \theta ~ (r_{i\pm1/2}-r^*_{i\pm1/2})$$ $$\delta Vol^\theta_{i,j \pm 1/2} =  r_{i} ~ \Delta r ~ (\theta_{j\pm1/2}-\theta^*_{j\pm1/2}).$$
Therefore we obtain:
\begin{equation}\label{}
%\left\{
\begin{array}{ll}
 \displaystyle  \int_{Vol_{i,j} }  f^{n+1}_{i,j}  ~ r dr d\theta =
 \displaystyle  \int_{Vol^n_{i,j} }  f(r,\theta,t^{n}) ~ rdr d\theta
\end{array}
\end{equation}
with $$ Vol^n_{i,j}= Vol_{i,j}-\delta Vol^r_{i+1/2}+\delta Vol^r_{i-1/2}-\delta Vol^\theta_{j+1/2}+\delta Vol^\theta_{j-1/2}. $$
We here recover a discrete mass conservation formulation. To obtain $Vol_{i,j}=Vol^n_{i,j}$, and thus preserve a constant function, it yields:
$$\delta Vol^r_{i+1/2}-\delta Vol^r_{i-1/2}+\delta Vol^\theta_{j+1/2}-\delta Vol^\theta_{j-1/2}=0.$$
Using $\delta Vol^k$ definitions,
we thus obtain a discrete divergence formulation in polar coordinates $\nabla^h \cdot a=0$ to be nullified:
\begin{equation}\label{div_disc_pol}
\begin{array}{cl}
\dfrac{1}{r_i}  \dfrac{  r_{i+1/2}  a_r(r_{i+1/2},\theta_j)- r_{i-1/2}  a_r(r_{i-1/2},\theta_j)}{\Delta r}+\\\dfrac{1}{r_i} \dfrac{r_i a_\theta(r_i,\theta_{j+1/2})- r_i a_\theta(r_i,\theta_{j-1/2})}{\Delta \theta}=0.
\end{array}
\end{equation}
with the following first order definition of the characteristic curves feet computation:
$$ r_{i\pm1/2}-r^*_{i\pm1/2}=\Delta t  ~a_r(r_{i\pm1/2},\theta_j) ~~ \mbox{and} ~~ \theta_{j\pm1/2}-\theta^*_{j\pm1/2}=\Delta t ~a_\theta(r_i,\theta_{j\pm1/2}).$$

\begin{figure} [!htbp]
\begin{center}
\includegraphics[height=4.5cm]{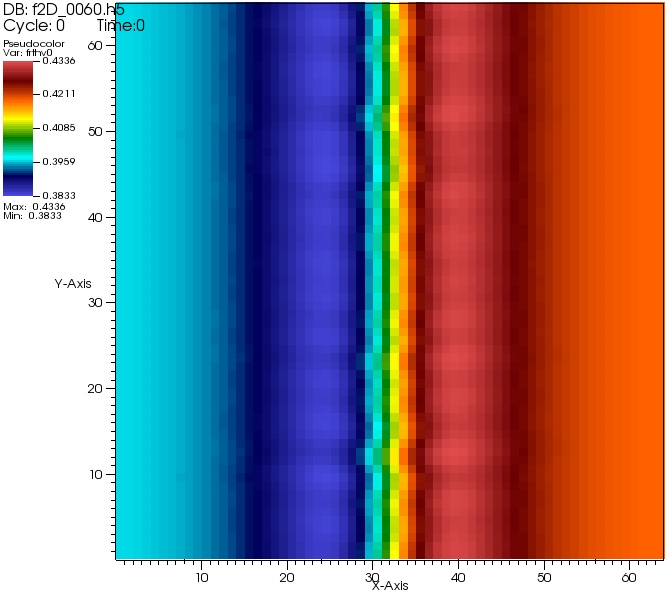}
\includegraphics[height=4.5cm]{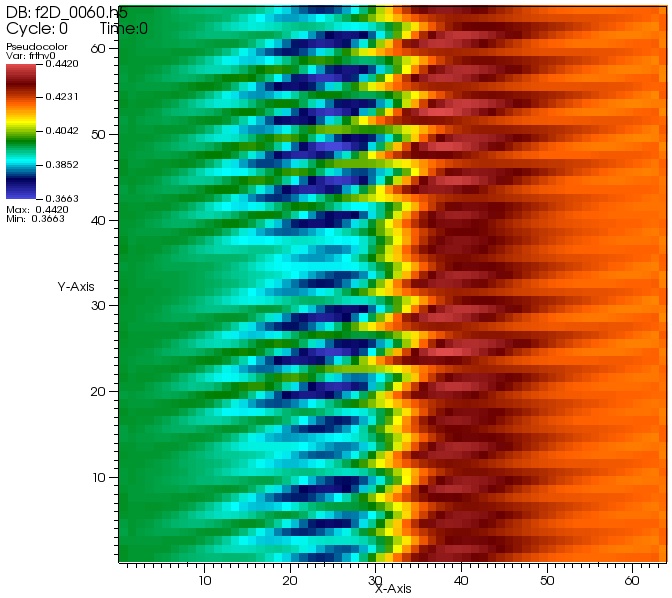}
\caption{\label{nodiv} Result at time $t=60$ with (left) the advection field computed in a way  (see section \ref{drift}) that satisfy the discrete divergence condition \eqref{div_disc_pol} and with (right) the advection field computed with cubic splines, which do not satisfy this condition \eqref{div_disc_pol}. Respecting condition \eqref{div_disc_pol} for the advection field not only leads to a better respect of the maximum principle, it is actually necessary to ensure the stability of the scheme. The result in figure \ref{nodiv} diverges from realistic physics. }
\end{center}
\end{figure}

As a conclusion, we have presented a general methodology for any coordinate system to compute the associated discrete divergence free condition, by using the approximation of cell edges by straight lines and by only considering the normal to cell faces motion of the volume as it has to be when invoking the {\it Green} formula in this Finite Volume framework.
The discrete divergence formulation  \eqref{div_disc_pol} is independent of the time integration method. It is a discrete consistent relation for the advection field of the form $\nabla^h \cdot a=0$. It should be satisfied to get the conservation condition on volumes $Vol_i^n=Vol_i^{n+1}$ in the phase space, which is necessary to obtain a maximum principle for the PSM scheme or actually for any Finite Volume scheme. This condition is  also necessary when using the Semi-Lagrangian PSM scheme with directional splitting as described in section \ref{section_PSM} as well as when using the Finite Volume form described in section \ref{PSM_FV}.  In Fig. \ref{nodiv}, we compare the results of a 4D drift-kinetic benchmark (see section \ref{drift} for details) obtained with the Semi-Lagrangian PSM scheme \ref{section_PSM} with an advection field computed: first in such a way the discrete divergence free condition \eqref{div_disc_pol} is satisfied and second with an advection field computed by cubic spline interpolation without satisfying this condition.
%#####################################################
\section{Use of the PSM scheme in a 4D drift-kinetic code}\label{use}

\subsection{Drift-kinetic model}\label{drift4d}
This work follows those of {\it Grandgirard et al} in the GYSELA code, see \cite{gysela} and \cite{gysela2}.
The geometrical assumptions of this model for ion plasma turbulence are a cylindrical geometry with $4D$ coordinates $(r,\theta,z,v_{\|})$ and a constant magnetic field $B=B_z ~e_z$,  where $e_z$ is the unit vector in $z$ direction. The model is the 4D Drift-Kinetic equations  described in{\it  Grandgirard et al} \cite{gysela2}:
\begin{equation}
\begin{array}{lll}
\dfrac{dr}{dt}=v_{GC_r}; ~r \dfrac{d \theta}{dt}=v_{GC_\theta}; ~ \dfrac{d z}{dt}=v_\|; ~\dfrac{d v_\|}{dt}=\dfrac{q_i}{m_i} E_z
\end{array}
\end{equation}
with $v_{GC}=(E \times B)/B^2$ and $E=-\nabla \Phi$ with $\Phi$ the electric potential.\\
The 4D Vlasov equation governing this system, where the ion distribution function is $f(r,\theta,z,v_\|,t)$, is the following:
\begin{equation}\label{Vlasov2}
\begin{array}{lll}
\partial_t f+v_{GC_r} \partial_r f + v_{GC_\theta} \partial_\theta f  + v_\| \partial_z f  +\dfrac{q_i}{m_i} E_z \partial_{v_\|} f =0.
\end{array}
\end{equation}
This equation is coupled with a quasi-neutrality equation for the electric potential $\Phi(r,\theta,z)$ that reads:
\begin{equation}\label{quasinetral}
\begin{array}{lll}
-\nabla_\perp \Phi \cdot \left( \dfrac{n_0(r)}{B~ \Omega_0} \nabla \Phi \right) + \dfrac{e~n_0(r)}{T_e(r)} (\Phi-<\Phi>_{\theta,z})=n_i-n_0
\end{array}
\end{equation}
with $n_i=\ds \int_{v_\|} f(r,\theta,z,v_\|) d v_\|$ and constant in time physical parameters $n_0$, $ \Omega_0$, $T_e$ and $e$.
Let us notice that the 4D velocity field $a=(v_{GC_r},v_{GC_\theta},v_\|, $ $q/m_i ~E_z)^t$ is divergence free:
\begin{equation}\label{diva}
\begin{array}{lll}
\nabla \cdot a = \dfrac{1}{r} \partial_r  (r~v_{GC_r})+  \dfrac{1}{r} \partial_\theta  (v_{GC_\theta})+\partial_z  v_\| + \partial_{v_\|} (q/m_i ~E_z)=0
\end{array}
\end{equation}
because of variable independence $\partial_{v_\|}  E_z=\partial_{v_\|} (\partial_z \Phi(r,\theta,z))=0$ and $\partial_z  v_\|=0$ and we have $v_{GC}=(E \times B)/B^2$, with  $E=-\nabla \Phi$ and $B=B_z ~e_z$, thus
\begin{equation}\label{def_vGC}
v_{GC_r}=\dfrac{1}{B_z} \left( -\dfrac{1}{r} \partial_\theta \Phi \right) ~~ \mbox{and} ~~ v_{GC_\theta}=\dfrac{1}{B_z} \left( \partial_r \Phi \right)
\end{equation}
and
\begin{equation}\label{div_rtheta}
\begin{array}{lll}
\nabla_{r, \theta} \cdot a &= \dfrac{1}{r}  \partial_r  (r~v_{GC_r})+ \dfrac{1}{r}  \partial_\theta  (v_{GC_\theta}) \\
~ & = \dfrac{1}{r~B_z}  \left( \partial_r  \left( r~ (-1/r) \partial_\theta \Phi \right) +  \partial_\theta \left( \partial_r \Phi \right) \right)=0.
\end{array}
\end{equation}
Therefore, one can write an equivalent conservative equation to the preceding Vlasov equation \eqref{Vlasov2}:
\begin{equation}\label{Vlasovc}
\begin{array}{lll}
\partial_t f+\partial_r (v_{GC_r} ~f)+  \partial_\theta (v_{GC_\theta} ~f)  +  \partial_z (v_\| ~f)  + \partial_{v_\|} \left( \dfrac{q_i}{m_i} E_z ~f \right) =0
\end{array}
\end{equation}

%-------------------------------------------------------------------------------------------------------
\subsection{Computation of a divergence free velocity field at the discrete level}
%-------------------------------------------------------------------------------------------------------
We have obtained a discrete form of the velocity field divergence to nullify (\ref{div_disc_pol}), as a necessary condition to obtain a numerical solution with a maximum principle. We saw in  (\ref{diva}) that $\nabla \cdot a =0$ is satisfied equivalently  if $\nabla_{r \theta} \cdot a =0$ (\ref{div_rtheta}) is satisfied and this is still true at the discrete level (independence of variables). Therefore, the velocity field should nullify the discrete polar divergence (\ref{div_disc_pol}):
\begin{equation}\label{divrt_disc}
\begin{array}{cl}
\dfrac{1}{r_i}  \dfrac{ r_{i+1/2}  a_r(r_{i+1/2},\theta_j)- r_{i-1/2}  a_r(r_{i-1/2},\theta_j)}{\Delta r}+\\\dfrac{1}{r_i} \dfrac{r_i a_\theta(r_i,\theta_{j+1/2})- r_i a_\theta(r_i,\theta_{j-1/2})}{\Delta \theta}=0
\end{array}
\end{equation}
with
$$a_r=dr / dt=v_{GC_r}=\dfrac{-1}{r~B_z} \partial_\theta \Phi ~~ \mbox{and}  ~~
 a_\theta=d \theta /dt=v_{GC_\theta}/r=\dfrac{1}{r~B_z}  \partial_r \Phi,$$
 using definitions given in (\ref{def_vGC}).
 \begin{prop}
 Let us define the electric potential at the nodes of the mesh $\Phi_{i+1/2,j+1/2}$, whatever the way it is computed.
Let us set the following natural finite difference approximation for the velocity field:
\begin{equation}\label{discvel}
\begin{array}{l}
a_r(r_{ i+1/2 },\theta_j)= \dfrac{ -1 }{ r_{i+1/2}~B_z} \dfrac{ \Phi_{i+1/2,j+1/2} - \Phi_{i+1/2,j-1/2} } {\Delta \theta}  \\
a_\theta(r_i,\theta_{j+1/2})=\dfrac{ 1 }{ r_i~B_z} \dfrac{ \Phi_{i+1/2,j+1/2} - \Phi_{i-1/2,j+1/2} } {\Delta r} .
\end{array}
\end{equation}
With this approximated velocity field, the approximation of $\nabla_{r \theta} \cdot a =0$  (\ref{divrt_disc}) is satisfied.
 \end{prop}
 The proof is easy, we just have to put the velocity field \eqref{discvel}  in (\ref{divrt_disc}) to see that all terms annulate each others.
 \begin{rem}
 Notice that the electric potential $\Phi$ should be computed at nodes $(i\pm1/2,j\pm1/2)$ of the mesh to obtain velocities at the center of cell faces $(i\pm1/2,j)$ and  $(i,j\pm1/2)$.  It is well adapted to the PSM schemes, where the displacement should be calculated at cell faces.  \end{rem}

%#####################################################
\subsection{Numerical tests}

%======================================================================
\subsubsection{Drift-kinetic 4D model, PSM schemes  comparison} \label{drift}
In this section, we will compare the numerical methods on a 4D drift-kinetic benchmark, following the paper {\it Grangirard et al} \cite{gysela2}. The model is described in section \ref{drift4d}. We will compute the growth of a 4D unstable  turbulent mode. The benchmark consists of exciting the plasma mode $(m,n)$, with $m$ the poloidal mode ($\theta$) and $n$ the toroidal mode ($z$).  The initial distribution function is the sum of an equilibrium and a perturbation distribution function $f=f_{eq}+\delta f$. The equilibrium distribution function has the following form:
\begin{equation}
f_{eq}(r,v_\|)=\dfrac{n_0(r)}{(2 \pi T_i(r) / m_i)^{1/2}} \exp \left( - \dfrac{m_i v_\|^2}{2 T_i(r)} \right)
\end{equation}
and the perturbation $\delta f$
\begin{equation}
\delta f(r,\theta,z,v_\|)=f_{eq}(r,v_\|) ~g(r)~ h(v_\|)~ \delta p(\theta,z)
\end{equation}
with $g(r)$ and $h(v_\|)$ two exponential functions and
$$ \delta p(\theta,z) =  \epsilon \cos \left(  \dfrac{2 \pi n}{L_z} z + m \theta  \right)$$
with $L_z$ the length of the domain in $z$ direction, $m_i$, $T_i(r)$, $n_0(r)$ physical constant profiles, see \cite{gysela2} for details.
We have set here $m=16$ and $n=8$.

%We run different tests comparing results obtained with the BSL scheme and the PSM scheme as described in this paper.
%\begin{figure}[!htbp]
%\begin{center}
%\includegraphics[height=6.cm]{IMAGES/noscale_100_BSL_128x256x128x64.jpg}
%\includegraphics[height=6.cm]{IMAGES/noscale_100_PSM_128x256x128x64.jpg}
%\caption{\label{} Simulation of 128x256x128x64 cells, comparison  BSL (left) and PSM (right), with color scales based on real computed data for each picture rather than set to the same color scale for all of them, time=2000.}
%\end{center}
%These pictures show that the BSL scheme and the PSM scheme do not satisfy a maximum principle, thus "singular" points may appear with values outside initial data bounds, here $[0.3873 ~; ~0.4415 ]$. \end{figure}

\subsubsection{Algorithm}
At the beginning of the time step, the distribution function $f(x,v_\|,t^n)$ is known at time $t^n$, with $x=(r,\theta,z)$. The time step $\Delta t=t^{n+1}-t^n$ is computed at each step with the CFL like condition:
$$\Delta t =  \displaystyle \min_{d=r,\theta,z,v_\|} \left( CFL_d ~\dfrac{\Delta x_d}{ \displaystyle \max_{x_d} (a^n_d(x_d))} \right),$$
with the coefficients $CFL_r=CFL_\theta=0.5$, because the flow is highly non-linear in $(r,\theta)$ planes thus characteristics should not cross each others during one time step, and $CFL_z=CFL_{v_\|}=8$ because it is linear advection in direction $z$ and $v_\|$ so characteristics can not cross each others and then we allow a maximum displacement of $8$ cells. Actually excluding the linear phase, the most restrictive directions for the time step are $r$ and $\theta$, in such a way this last value (8) has a minor importance compare to the leading parameters  $CFL_r=CFL_\theta=0.5$. \\
The operator splitting between the quasi-neutral equation and the Vlasov transport equation is made second order using a Predictor-Corrector scheme in time:
\begin{enumerate}
\item Time step $\Delta t$ computation.
\item Quasi-neutral equation \eqref{quasinetral} solving at $t^n$ using the distribution function $f^n$ (actually the density) to obtain the electric potential $\Phi^{n}(x)$ at time $t^{n+1/2}$. The advection field $a(x,t^{n})$ is computed with $\Phi^{n}(x)$  according to equation \eqref{diva} and using formula \eqref{discvel} .
\item 4D Vlasov equation solving at $t^n$ with time step $\Delta t / 2$ to obtain the distribution function $f^{n+1/2}(x)$ at time $t^{n+1/2}$ using the advection field $a(x,t^{n})$.
\item Quasi-neutral equation  \eqref{quasinetral} solving at $t^{n+1/2}$  using the distribution function $f^{n+1/2}$ (actually the density) to obtain the electric potential $\Phi^{n+1/2}(x)$ at time $t^{n+1/2}$. The advection field $a(x,t^{n+1/2})$ is computed with $\Phi^{n+1/2}(x)$  according to equation \eqref{diva} and using formula \eqref{discvel} .
\item 4D Vlasov equation solving at $t^n$ with time step $\Delta t$ to obtain the distribution function $f^{n+1}(x)$ at time $t^{n+1}$ using the advection field $a(x,t^{n+1/2})$.
\end{enumerate}

In the two following paragraphs, we describe the schemes for the 4D Vlasov equation \eqref{Vlasovc} solving of the algorithm with $\Delta t^*=\Delta t/2$ in the prediction step and  $\Delta t^*=\Delta t$ in the correction step.

%-------------------------------------------------------------------------------------------------------
\paragraph{4D Semi-Lagrangian PSM sheme with directional splitting}\label{algosl}
%-------------------------------------------------------------------------------------------------------
\begin{itemize}
\item PSM 1D advection of $f(x,v_\|,t^n)$ in direction $v_\|$ with velocity $a^n_{v_\|}$ and time step $\Delta t^*/2$ to obtain $f(x,v_\|,t^{v_\| /2})$.
\item PSM 1D advection of $f(x,v_\|,t^{v_\| /2})$ in direction $z$ with velocity $a^n_{z}$ and time step $\Delta t^*/2$ to obtain $f(x,v_\|,t^{z /2})$.
\item PSM 1D advection of $f(x,v_\|,t^{z /2})$ in direction $\theta$ with velocity $a^n_{\theta}$ and time step $\Delta t^*/2$ to obtain $f(x,v_\|,t^{\theta /2})$.
\item PSM 1D advection of $f(x,v_\|,t^{\theta /2})$ in direction $r$ with velocity $a^n_{r}$ and time step $\Delta t^*$ to obtain $f(x,v_\|,t^{r})$.
\item PSM 1D advection of $f(x,v_\|,t^r)$ in direction $\theta$ with velocity $a^n_{\theta}$ and time step $\Delta t^*/2$ to obtain $f(x,v_\|,t^{\theta})$.
\item PSM 1D advection of $f(x,v_\|,t^\theta)$ in direction $z$ with velocity $a^n_{z}$ and time step $\Delta t^*/2$ to obtain $f(x,v_\|,t^{z})$.
\item PSM 1D advection of $f(x,v_\|,t^z)$ in direction $v_\|$ with velocity $a^n_{v_\|}$ and time step $\Delta t^*/2$ to obtain $f(x,v_\|,t^{v_\|})=f(x,v_\|,t^{n+1})$.
 \end{itemize}
Each PSM 1D advection is achieved using the standard 1D semi-Lagrangian PSM scheme as described in section \ref{section_PSM}. The directional splitting is second order by using a Strang like decomposition.  Since we use here a directional splitting and a second order scheme in time for the computation of the characteristic curves, the volumes are not strictly conserved in the phase space, because the scheme does not satisfy  the discrete  divergence free condition \eqref{div_disc_pol}.

%-------------------------------------------------------------------------------------------------------
\paragraph{4D Finite Volume form of the PSM scheme:}\label{algofv}
%-------------------------------------------------------------------------------------------------------
%At the beginning of the time step, the distribution function $f(x,v_\|,t^n)$ is known at time $t^n$, with $x=(r,\theta,z)$. The time step is $\Delta t^*=t^{n+1}-t^n$.
 \begin{itemize}
%\item The quasi-neutral equation is resolved using the distribution function $f(x,v_\|,t^n)$ to obtain the electric potential $\Phi^n(x)$ at time $t^n$.
%\item The advection field $a(x,t^n)=(a^n_{r}(x), a^n_{\theta}(x), a^n_{z}(x), a^n_{v_\|}(x))^t $ is computed with $\Phi^n(x)$ according to equation \eqref{diva} and using formula \eqref{discvel} .
\item PSM 1D advection of $f(x,v_\|,t^n)$ in direction $v_\|$ with velocity $a^n_{v_\|}$ and time step $\Delta t^*/2$ to obtain $f(x,v_\|,t^{v_\| /2})$.
\item PSM 1D advection of $f(x,v_\|,t^{v_\| /2})$ in direction $z$ with velocity $a^n_{z}$ and time step $\Delta t^*/2$ to obtain $f(x,v_\|,t^{z /2})$.
\item PSM 2D advection of $f(x,v_\|,t^{z /2})$ in each plane $(r,\theta)$ with velocities $(a^n_{r},a^n_{\theta})$ and time step $\Delta t^*$ to obtain $f(x,v_\|,t^{r,\theta })$.
\item PSM 1D advection of $f(x,v_\|,t^{r,\theta})$ in direction $z$ with velocity $a^n_{z}$ and time step $\Delta t^*/2$ to obtain $f(x,v_\|,t^{z})$.
\item PSM 1D advection of $f(x,v_\|,t^z)$ in direction $v_\|$ with velocity $a^n_{v_\|}$ and time step $\Delta t^*/2$ to obtain $f(x,v_\|,t^{v_\|})=f(x,v_\|,t^{n+1})$.
 \end{itemize}
Each PSM 1D advection is achieved using the standard 1D semi-Lagrangian scheme as described in section \ref{section_PSM}. The PSM 2D advection in $(r,\theta)$ is achieved with the Finite Volume form as described in section \ref{PSM_FV}. Since we use here the scheme  \ref{PSM_FV}, the volumes are strictly conserved in the phase space, because the scheme does satisfy  the discrete divergence free condition \eqref{div_disc_pol}. Even if we use the semi-Lagrangian PSM 1D advection in $z$ and $v_\|$ directions, the property is kept because the velocity is constant in these directions.

%-------------------------------------------------------------------------------------------------------
\subsection{Results}
\begin{figure}[!htbp]
\begin{center}
\vspace{-2cm}
{\bf PSM Directional Splitting}~~~~~~~~~~~~~{\bf  PSM Finite Volume}\\
\includegraphics[height=9.cm]{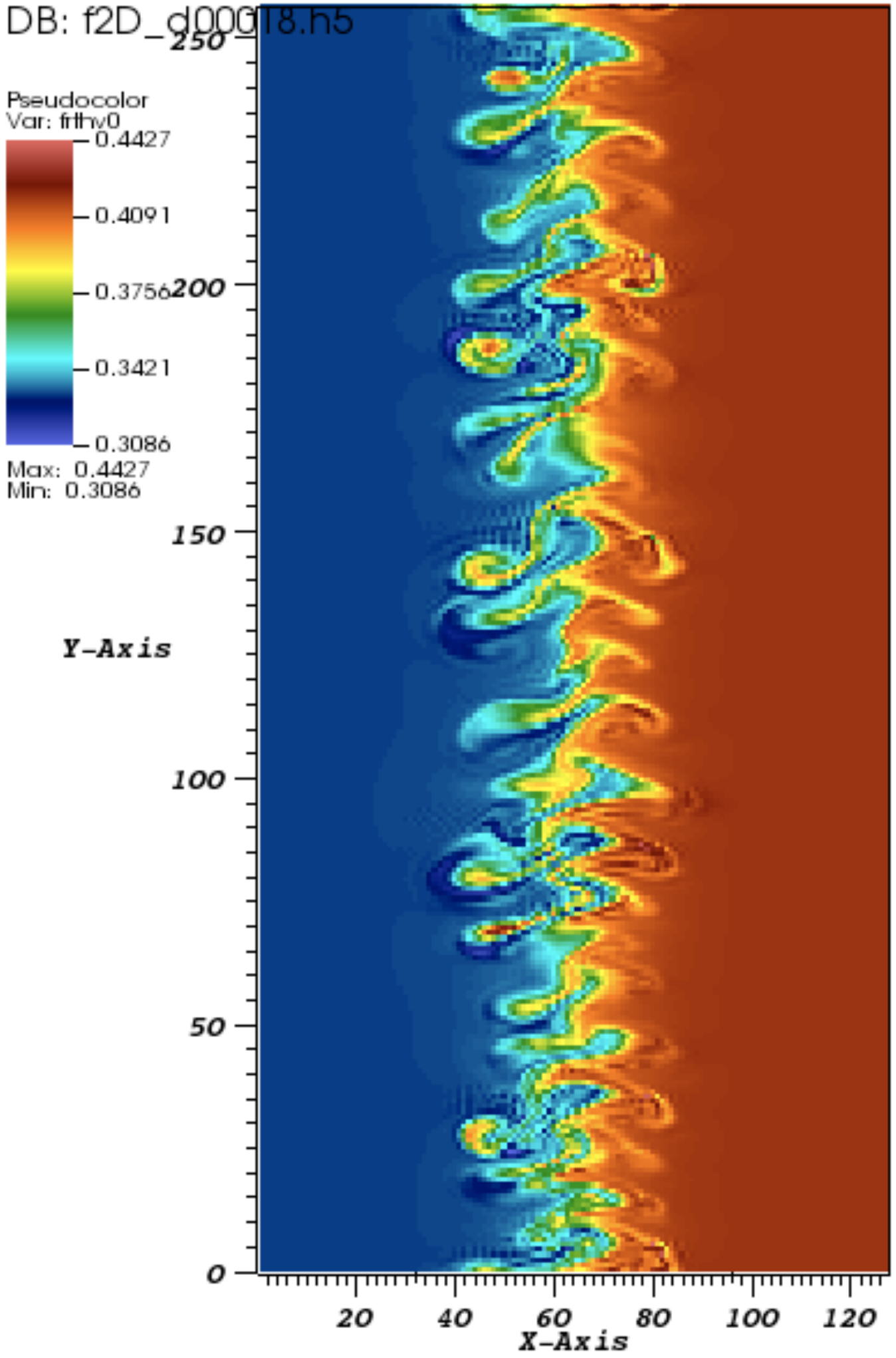}
\includegraphics[height=9.cm]{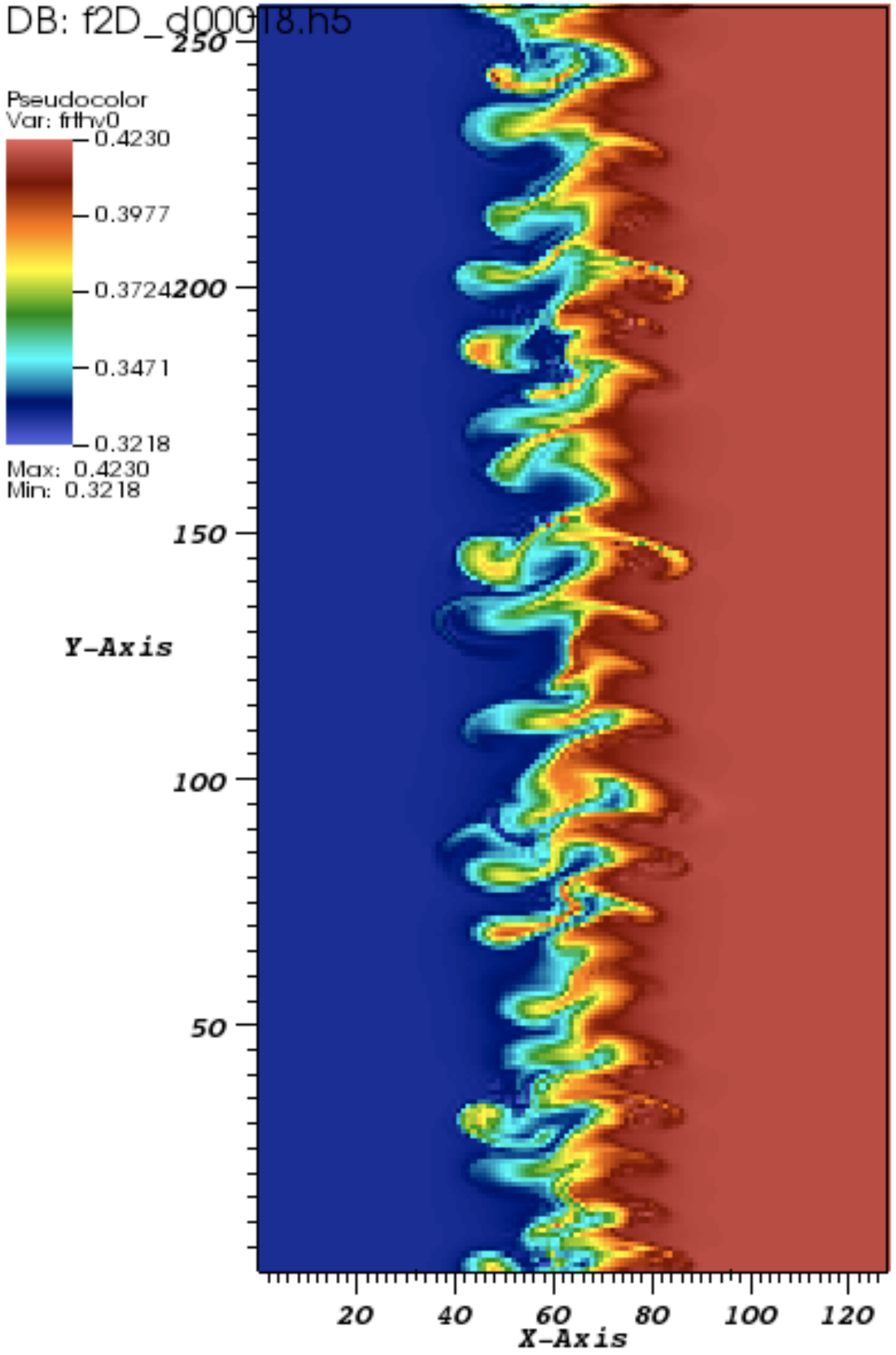}\\
{\bf SLS Directional Splitting}~~~~~~~~~~~~~ {\bf SLS Finite Volume}\\
\includegraphics[height=9.cm]{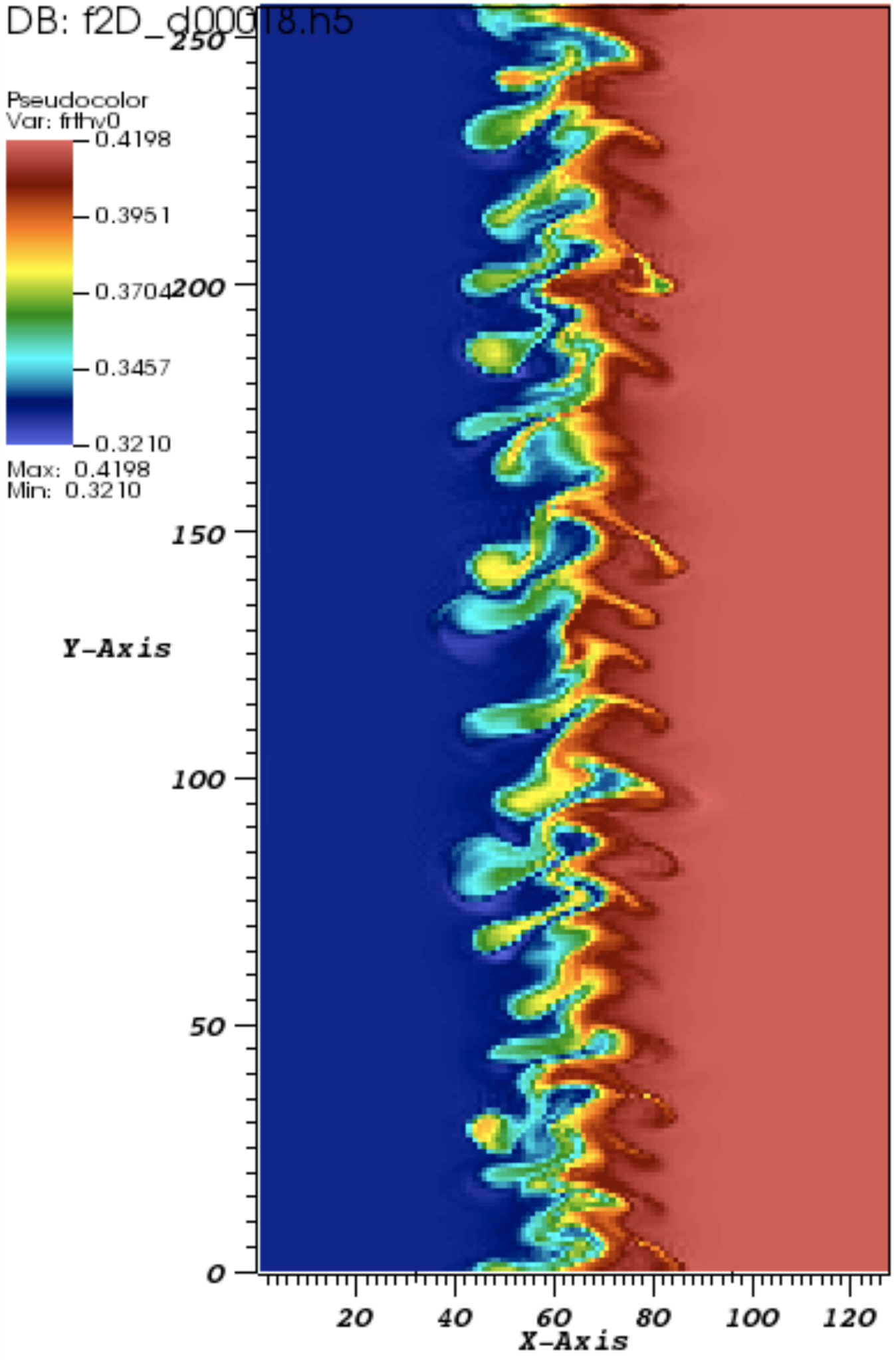}
\includegraphics[height=9.cm]{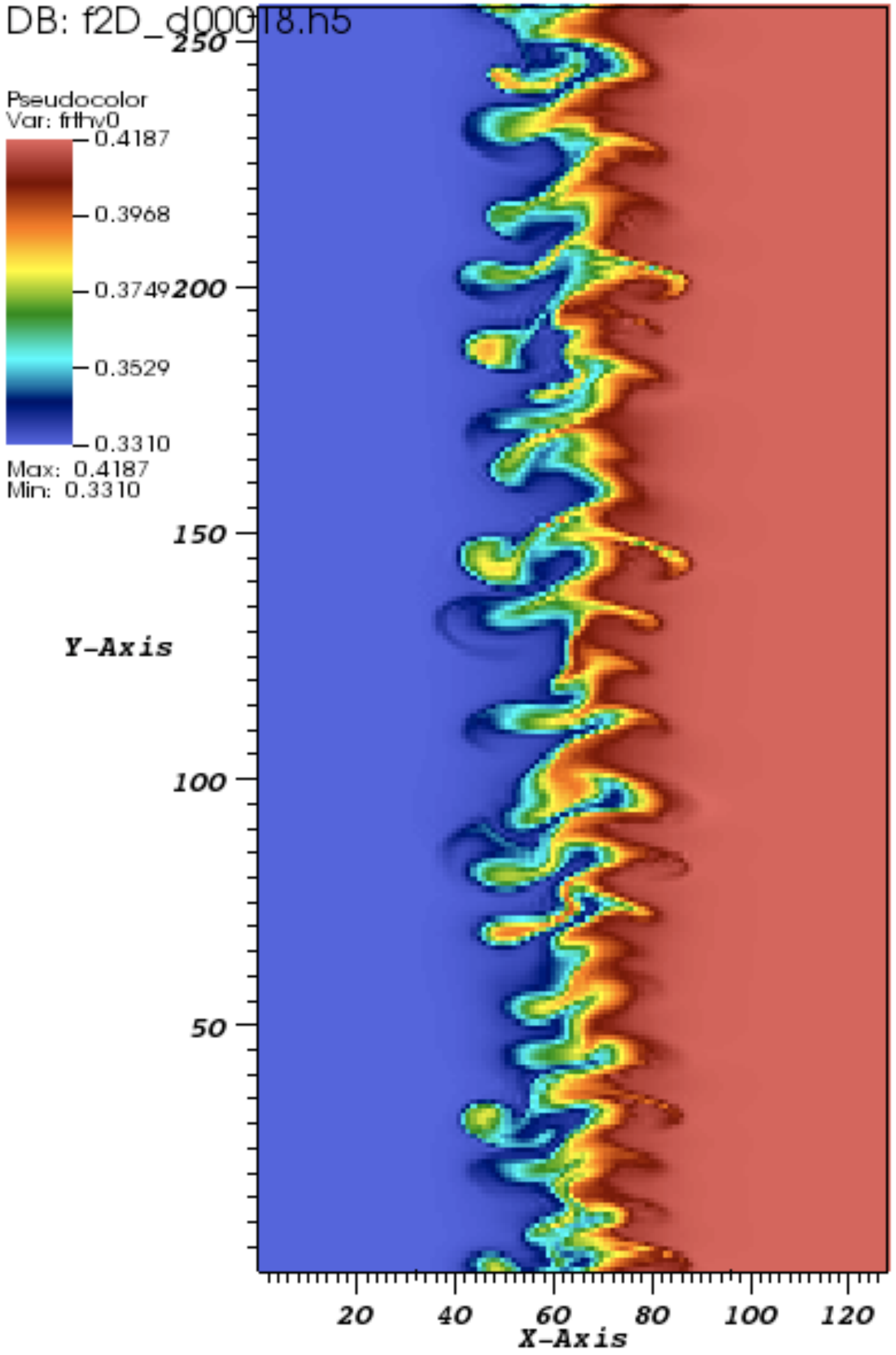}\\
\caption{\label{sol1800} Simulation with 128x256x32x16 cells  --- PSM Directional Splitting 1D (up-left) ---  PSM Finite Volume (up-right) --- SLS Directional Splitting 1D (down-left) --- SLS Finite Volume (down-right) ---  time =1800.}
\end{center}
\end{figure}

\begin{figure}[!htbp]
\begin{center}
\vspace{-2cm}
{\bf PSM Directional Splitting}~~~~~~~~~~~~~{\bf  PSM Finite Volume}\\
\includegraphics[height=9.cm]{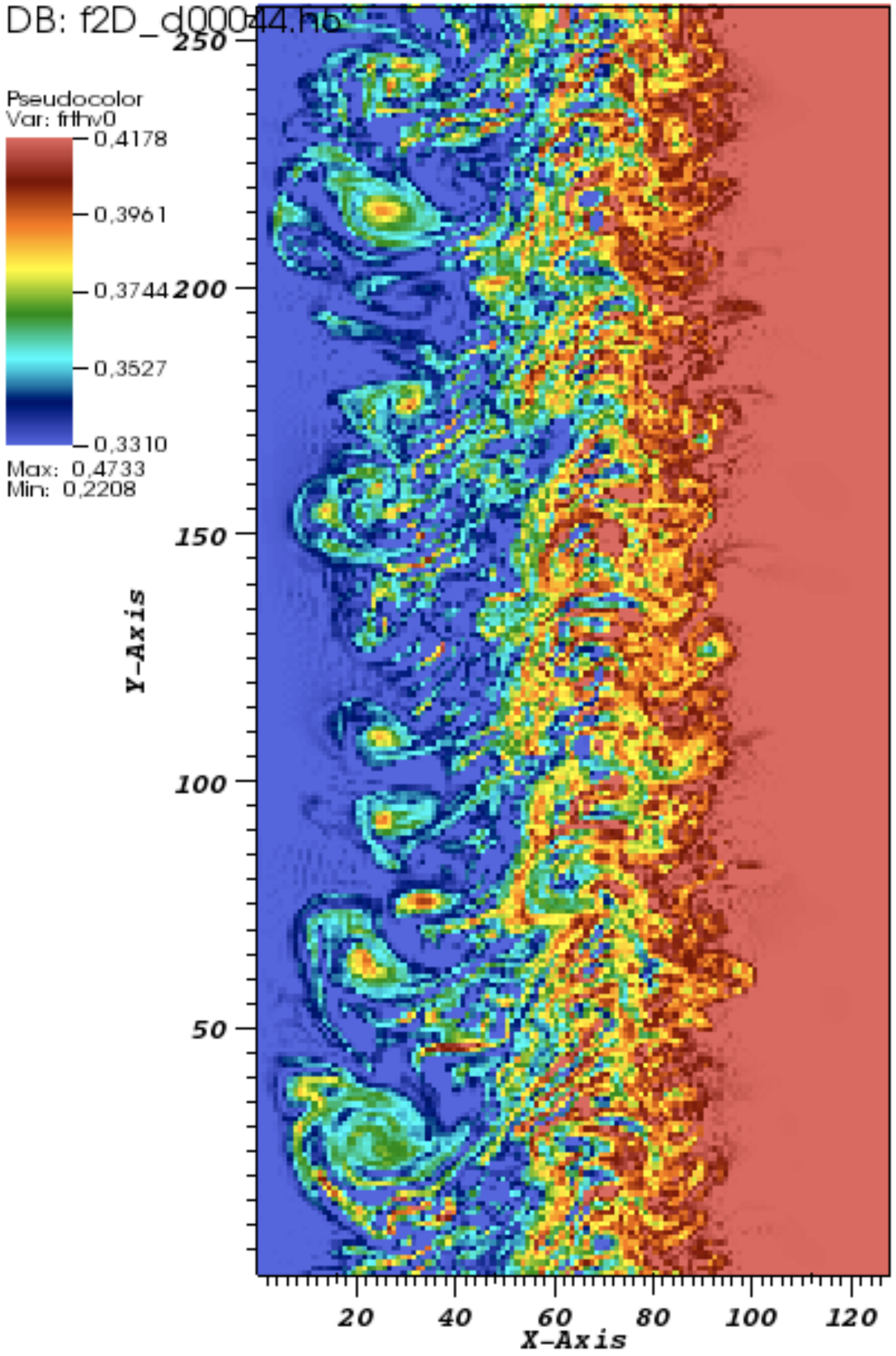}
\includegraphics[height=9.cm]{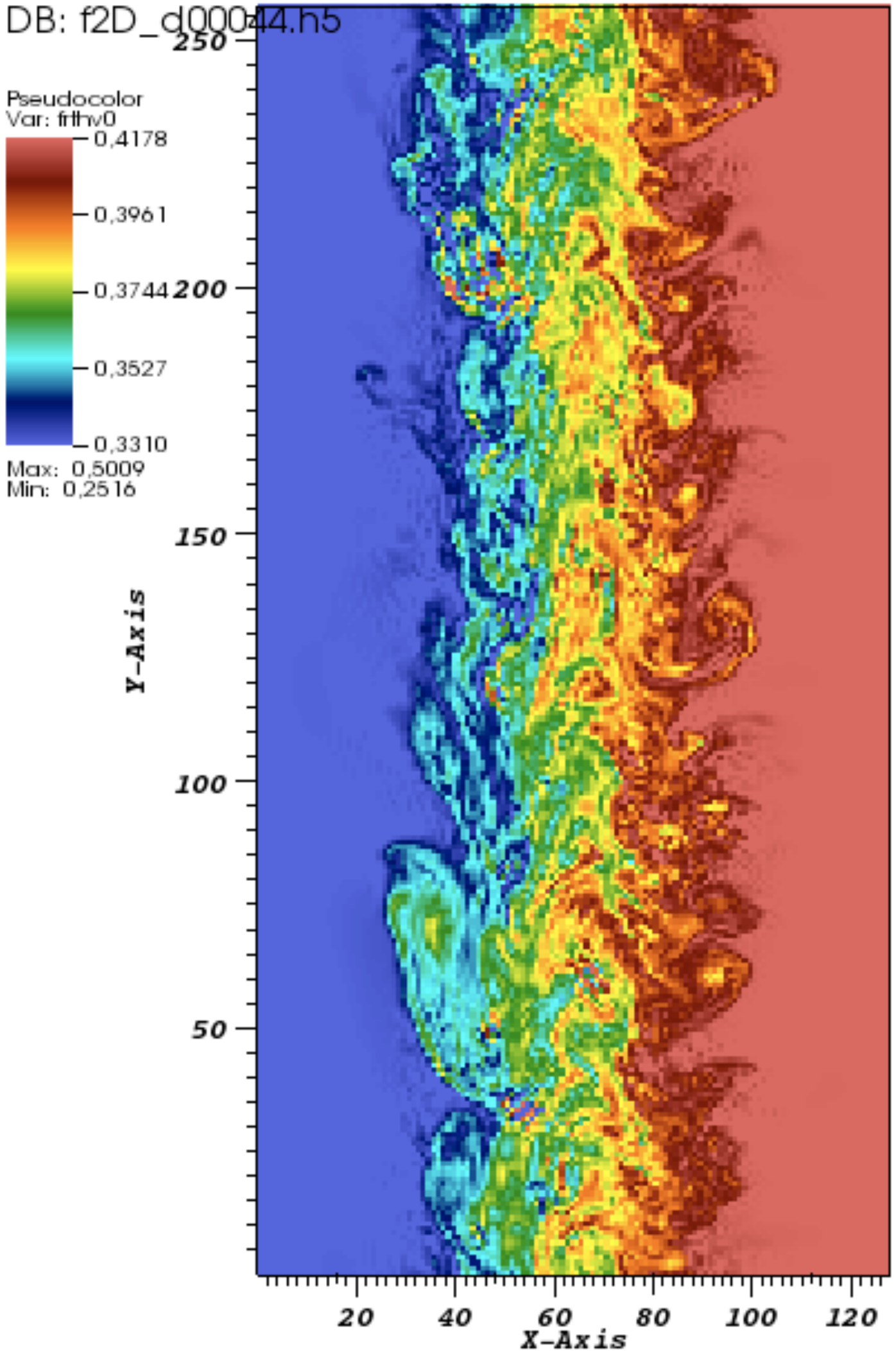}\\
{\bf SLS Directional Splitting}~~~~~~~~~~~~~ {\bf SLS Finite Volume}\\
\includegraphics[height=9.cm]{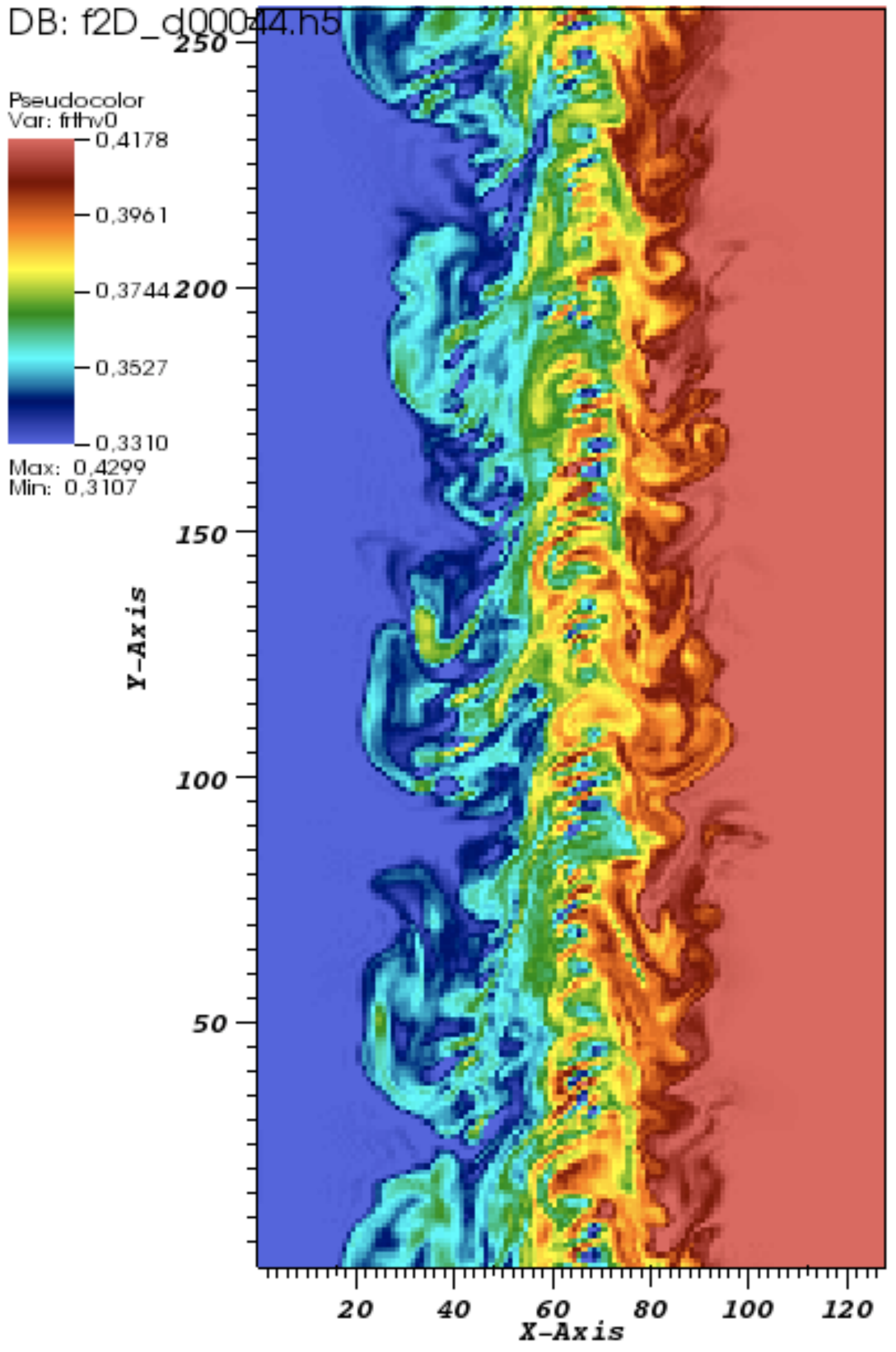}
\includegraphics[height=9.cm]{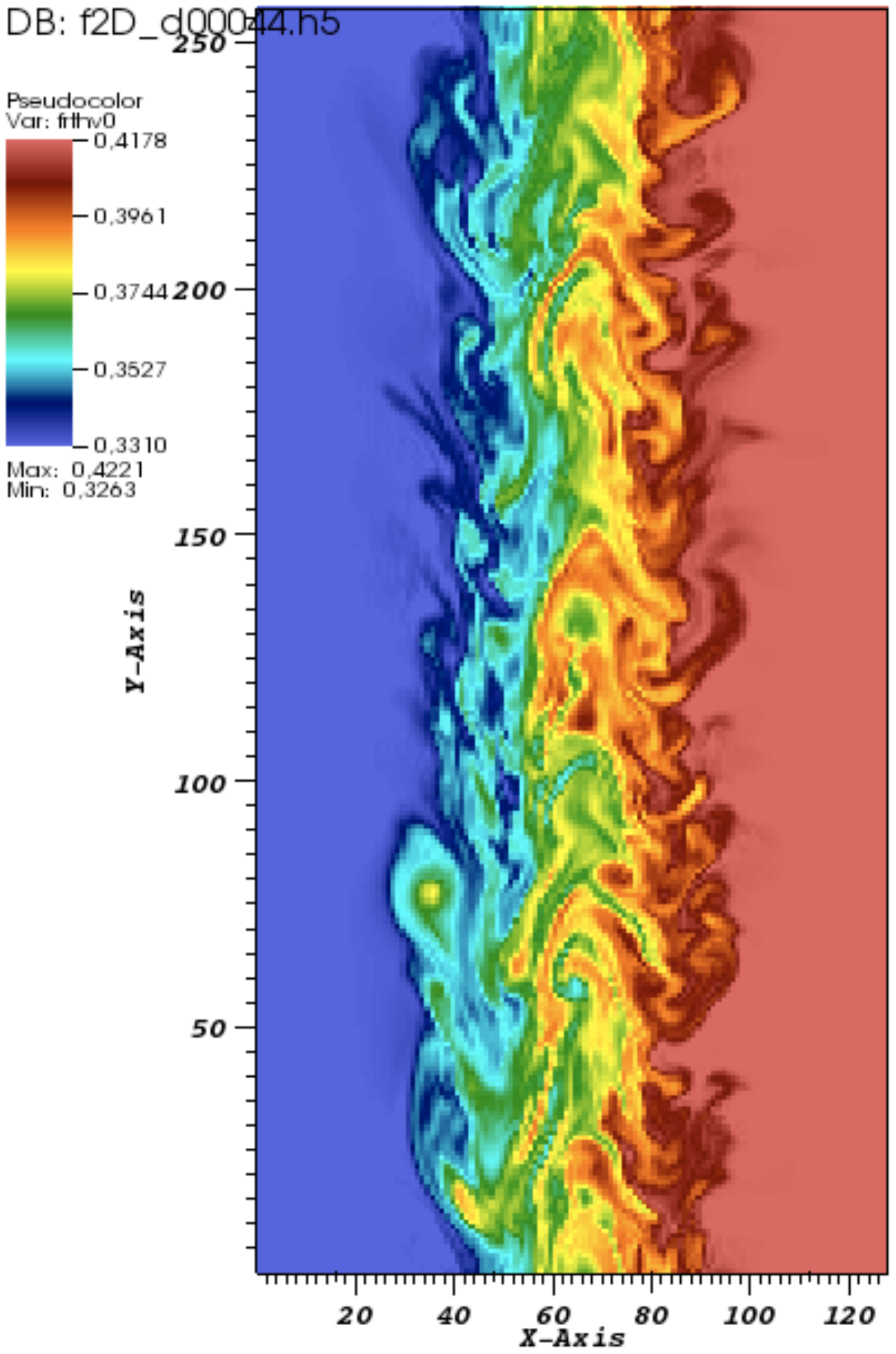}
\caption{\label{sol4400} Simulation with 128x256x32x16 cells  --- PSM Directional Splitting 1D (up-left) ---  PSM Finite Volume (up-right) --- SLS Directional Splitting 1D (down-left) --- SLS Finite Volume (down-right) ---  time =4400, with all color tables set to the minimum and maximum value at initial time.}
\end{center}
\end{figure}

\begin{figure}[!htbp]
\begin{center}
\vspace{-2cm}
{\bf BSL  time=1800}~~~~~~~~~~~~~~~~~~~~~~~~~{\bf  BSL time=4400}\\
\includegraphics[height=9.cm]{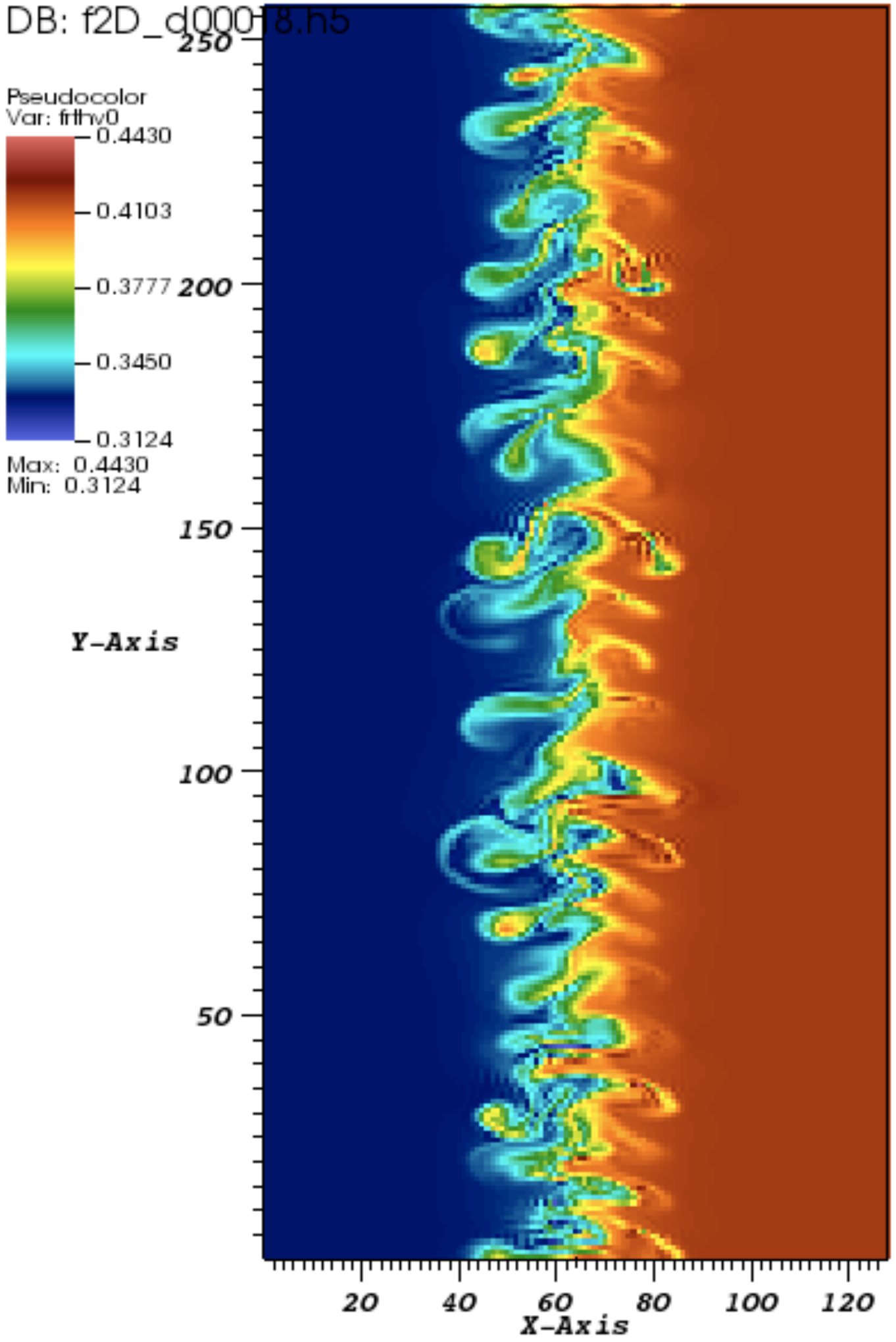}
\includegraphics[height=9.cm]{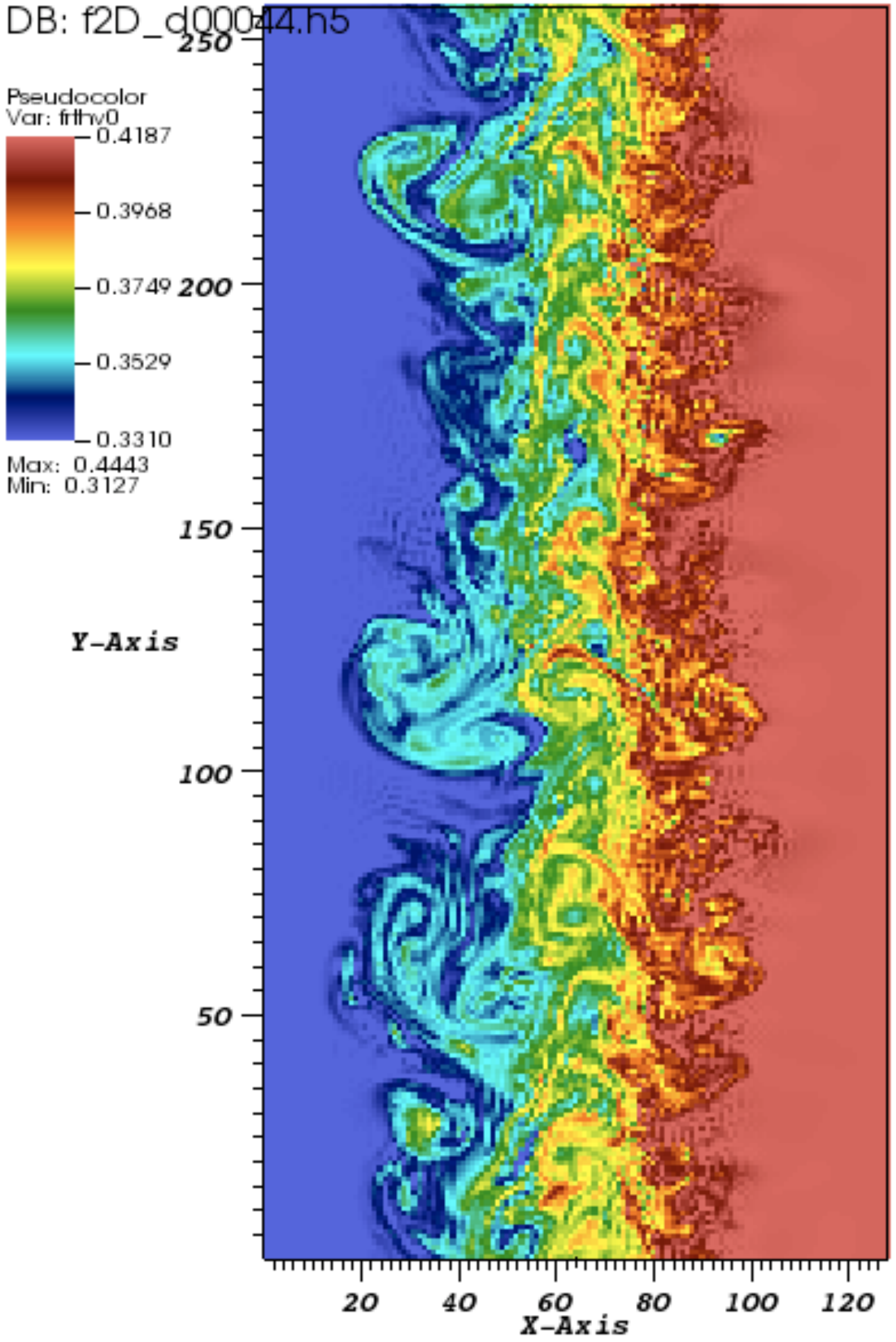}\\
\caption{\label{solBSL} Reference BSL Simulation with 128x256x32x16 cells  --- BSL at time=1800 with the real color table values  (left) ---  BSL at time=4400 with the color table set to the minimum and maximum value at initial time (left).}
\end{center}
\end{figure}

%-------------------------------------------------------------------------------------------------------
The mesh is $128 \times256 \times32 \times16$ cells in $r,\theta,z,v_\|$ directions. Boundary conditions are periodic for directions $\theta$ and $z$ and {\it Neumann} ($\partial f/ \partial n=0$) in $r$ and $v_\|$.
We first ran the reference test case with the non-conservative Backward Semi-Lagrangian (BSL) scheme  in section \ref{sec:BSLscheme} which is currently used in the GYSELA code.
Then we ran four test cases to show the influence of each numerical treatment: the standard conservative semi-Lagrangian PSM scheme in section \ref{algosl}  with 1D directional splitting (PSM Directional Splitting 1D) and the same with the SLS limiter (SLS Directional Splitting 1D),  the unsplit Finite Volume form of the PSM scheme in section \ref{algofv} (PSM Finite Volume) and the same with the SLS limiter (SLS Finite Volume).\\
The computed 4D distribution functions are pictured in Fig. \ref{sol1800} at time $t=1800$ and in Fig. \ref{sol4400} at time $t=4400$. We only present 2D slices $(r, \theta)$ of the distribution function at $v_\|=0$ and for a given value of $z=z_0$. In these figures $X$ stands for $r$ direction and $Y$ for the $\theta$ direction.
At initial time $t=0$, the minimum and maximum values of the distribution function in this slice are $(\min.=0.331, \max.=0.4187)$ and these values should be the same at any time of the computation if the maximum principle would be respected.
In Fig. \ref{sol1800}, we show pictures of each scheme result at time $t=1800$, which corresponds approximately to the beginning of the non-linear turbulent phase saturation. Small structures are appearing and interact with each others. All results are still close qualitatively. However, we already see oscillations in the solution obtained with PSM DS (Directional Splitting), where the minimum and maximum values $(\min.=0.3086, \max.=0.4427)$ are already quite different than the one at initial time. The PSM Finite Volume (PSM FV) form and the SLS DS better keep these extrema, but only SLS Finite Volume keep the extrema unchanged until time $t=1800$ with a really similar behaviour of the solution.
In Fig. \ref{sol4400}, we show pictures of each scheme result at time $t=4400$ when turbulence is well developed. We see that the standard PSM scheme creates a lot of unphysical oscillations (structures are reaching the boundaries in $r$) and may crash the computation. The PSM FV form and the SLS DS better keep the turbulence structures, but still oscillations are created. The SLS Finite Volume keep the extrema of the solution reasonably well $(\min.=0.3263, \max.=0.4221)$ (the SLS limiter does not provide a maximum principle) and the solution is smooth. We may say that the added diffusion with the limiter helps the scheme to diffuse subgrid structures without creating oscillations. The divergence free property of the Finite Volume scheme is important to cure to solution from instabilities that can be seen at $r$ values close the average value of $r$ (vertical line at the middle of pictures  in Fig. \ref{sol4400}) in the SLS Directional Splitting solution compare to the  SLS Finite Volume solution.
In Fig. \ref{solBSL}, we see in the reference BSL solution at time $t=1800$ spurious oscillations produced during the reconstruction step of the distribution function, which is the only possibility to break the maximum principle for the BSL scheme: here extrema are $(\min.=0.3124, \max.=0.4430)$ instead of values at initial time $(\min.=0.331, \max.=0.4187)$. At time $t=4400$, we see spurious oscillations as well, but the maximum principle is better satisfied than with the standard PSM DS  scheme in  Fig. \ref{sol4400}, because no conservation of volumes in the phase space  has to be satisfied, as it is explained in section \ref{sec:BSLscheme}.
%\newpage
%========================================================================
%   C O N C L U S I O N
\section{Conclusion and perspectives} \label{sec:conclusion}
% Conclusion
% Perspectives
The PSM scheme has been successfully integrated in the GYSELA code and has been tested on 4D Drift-Kinetic test cases. We had first experimentally stated and afterward explained in this paper that the PSM scheme can be unstable without taking care of a velocity field divergence free condition.  The numerical results show that the study of the volume evolution in the phase space is fruitful. Notice that this conservative scheme properly allows a directional splitting, in the semi-Lagrangian or in the Finite Volume form,  what is not the case with the BSL scheme. The Slope Limited Splines (SLS) limiter is efficient to cut off spurious oscillations of the standard PSM scheme by adding diffusion that helps eventually the scheme to manage small structures below the cell size.  Of course, the PSM scheme should be further validated as well as its integration in the GYSELA code using the gyrokinetic 5D model in toroidal geometry. In particular, the curvature of the mesh couple several directions  by the geometrical  {\it Jacobian} which makes the divergence free condition more complex, as well as the writing of the Quasi-Neutral solver and the Gyroaverage operator.

%========================================================================

%========================================================================
%  B I B L I O G R A P H Y
 
%========================================================================

\newpage
\tableofcontents

\end{document}